\documentclass[11pt,a4paper]{amsart}
\usepackage{amssymb,amscd}
\usepackage[all]{xy}

\numberwithin{equation}{section}
\newtheorem{theorem}{Theorem}[section]
\newtheorem{lemma}[theorem]{Lemma}
\newtheorem{prop}[theorem]{Proposition}
\newtheorem{corol}[theorem]{Corollary}
\newtheorem{conj}[theorem]{Conjecture}
\newtheorem{step}{Step}
\theoremstyle{definition}
\newtheorem{defin}[theorem]{Definition}

\theoremstyle{remark}
\newtheorem{remk}[theorem]{Remark}
\newtheorem*{erem}{Remark}

\def\labelenumi{\theenumi.}
\def\brn{\def\labelenumi{\theenumi)}}

\def\iff{if and only if }
\def\CM{Cohen--Mac\-au\-lay }
\def\VB{vector bundle}
\def\oc{one-to-one correspondence}

\def\RR{Riemann--Roch Theorem}
\def\dvr{discrete valuation ring}

\def\op{^\mathrm{op}}
\def\chr{\mathop\mathrm{char}\nolimits}
\def\ting{\mathop{\widetilde{\mathrm{Sing}}}\nolimits}
\def\rk{\mathop\mathrm{rk}\nolimits}
\def\he{\mathop\mathrm{ht}\nolimits}
\def\pic{\mathop\mathrm{Pic}\nolimits}
\def\hom{\mathop\mathrm{Hom}\nolimits}

\def\ext{\mathop\mathrm{Ext}\nolimits}
\def\Hom{\mathop{\mathcal H\!\mathit{om}}\nolimits}

\def\spec{\mathop\mathrm{Spec}\nolimits}
\def\sing{\mathop\mathrm{Sing}\nolimits}
\def\Vb{\mathop\mathrm{VB}\nolimits}
\def\Cm{\mathop\mathrm{CM}\nolimits}
\def\fsm{\mathop\mathrm{Fsm}\nolimits}
\def\pro#1{\mathrm{pro}\mbox{-}#1}
\def\md#1{#1\mbox{-}\mathrm{mod}}
\def\El{\mathop\mathrm{El}\nolimits}
\def\rad{\mathop\mathrm{rad}\nolimits}
\def\res{\mathop\mathrm{res}\nolimits}
\def\ann{\mathop\mathrm{Ann}\nolimits}
\def\GL{\mathop\mathrm{GL}\nolimits}
\def\diag{\mathop\mathrm{diag}\nolimits}

\def\poly{\mathop\mathrm{Poly}\nolimits}
\def\chr{\mathop\mathrm{char}\nolimits}
\def\mat{\mathop\mathrm{Mat}\nolimits}
\def\gdim{\mathop\mathrm{gl.dim}\nolimits}
\def\syz{\mathop\mathrm{syz}\nolimits}

\def\rap{\noindent}	\def\bul{\bullet}
\def\ti{\tilde}	
\def\vv{^\vee}	\def\vvv{^{\vee\vee}}
\def\vs{^\sharp}	\def\mps{\mapsto}
\def\lst#1#2{#1_1,#1_2,\dots,#1_{#2}}
\def\row#1#2{(#1_1,#1_2,\dots,#1_{#2})}
\def\set#1{\left\{\,#1\,\right\}}
\def\setsuch#1#2{\left\{\,#1\,|\,#2\,\right\}}
\def\gnr#1{\langle\,#1\,\rangle}

\def\mtr#1{\begin{pmatrix}#1\end{pmatrix}}
\def\larr{\longrightarrow}
	
\def\str{\stackrel}      \def\bsc{\bigsqcup}
\def\bap{\bigcap}        \def\bup{\bigcup}
\def\sb{\subset}         \def\sp{\supset}
      \def\sbe{\subseteq}

        \def\ul{\underline}
	\def\bop{\bigoplus}

\def\={\setminus}        \def\8{\infty}
\def\0{\emptyset}        \def\3{\ss}
\def\*{\otimes}          
\def\ti{\widetilde}      \def\iso{\simeq}
        \def\xx{\times}
\def\bop{\bigoplus}      
           \def\+{\oplus}
\def\mdl#1{#1\mbox{-}\mathrm{mod}}
\def\coh{\mathop\mathrm{Coh}\nolimits}
\def\bot{\boxtimes}		
\def\pro#1{#1\mbox{-}\mathrm{pro}}

\def\Mk{\Bbbk}

\def\la{\lambda}	\def\al{\alpha}
	\def\ga{\gamma}
\def\La{\Lambda}	\def\Ga{\Gamma}
\def\De{\Delta}		\def\Si{\Sigma}
\def\de{\delta}		\def\om{\omega}
\def\eps{\varepsilon}	\def\si{\sigma}
		\def\th{\theta}
\def\Om{\Omega}	\def\vi{\varphi}

\def\bA{\mathbf A}	\def\fB{\mathbf b}
\def\fD{\mathbf d}
\def\bF{\mathbf F}	\def\bP{\mathbf P}
\def\bV{\mathbf V}	\def\fC{\mathbf c}
\def\fU{\mathbf u}	\def\fV{\mathbf v}
\def\fE{\mathbf e}	\def\bR{\mathbf R}
\def\bB{\mathbf B}	\def\bK{\mathbf K}
\def\bL{\mathbf L}	
	
\def\fP{\mathbf p}	\def\bT{\mathbf T}

\def\mN{\mathbb N}	\def\mP{\mathbb P}
\def\mZ{\mathbb Z}		\def\mC{\mathbb C}
\def\mA{\mathbb A}	

\def\rA{\mathrm A}	\def\rH{\mathrm H}
\def\rD{\mathrm D}	\def\rE{\mathrm E}
\def\rT{\mathrm T}	\def\nH{\mathrm h}
	\def\rT{\mathrm T}
	
\def\kO{\mathcal O}	\def\tO{\tilde{\mathcal O}}
\def\kF{\mathcal F}	\def\kE{\mathcal E}
\def\kA{\mathcal A}	\def\kB{\mathcal B}
\def\kV{\mathcal V}	\def\kG{\mathcal G}
\def\kP{\mathcal P}	\def\kS{\mathcal S}	
\def\kJ{\mathcal J}	\def\kL{\mathcal L}
\def\kM{\mathcal M}\def\kN{\mathcal N}
\def\tG{\tilde{\mathcal G}}

\def\sA{\mathsf A}	\def\sB{\mathsf B}
\def\sU{\mathsf U}

\def\gM{\mathfrak m} 	\def\dE{\mathfrak E}
\def\dF{\mathfrak F}	\def\dM{\mathfrak M}
\def\dN{\mathfrak N}	\def\dI{\mathfrak I}
\def\gP{\mathfrak p}	\def\gQ{\mathfrak q}
\def\gX{\mathfrak x}
	
\def\tF{\tilde{\mathcal F}}
\def\oX{\breve X}	\def\tX{\tilde X}	
\def\tbF{\tilde{\mathbf F}}
\def\tV{\tilde{\mathcal V}}
\def\oS{\breve S}

\begin{document}

\title{Vector bundles and Cohen--Macaulay modules}
\author{Yuriy A. Drozd}
\address{Kiev Taras Shevchenko University,
 Department of Mechanics and Mathematics, 01033 Kiev, Ukraine, and
 Universit\"at Kaiserslautern, Fachbereich Mathematik, 67663 Kaiserslautern, Germany}
\email{unialg@ln.ua}
\thanks{Supported by the  DFG Schwerpunkt
 ``Globale Methoden in der komplexen Geometrie'' and by the CRDF Award UM\,2-2094}
\dedicatory{To the memory of Sheila Brenner}
\subjclass[2000]{Primary 13C14, 14H60; Secondary 13H10, 14J17}
\maketitle

 \tableofcontents
 \section*{Introduction}

  The aim of this survey is to present recent results on classification of vector bundles over projective curves and Cohen--Macaulay modules over
 surface singularities, mainly obtained by the author in collaboration with G.-M.\,Greuel and I.\,Kashuba \cite{dg3,dgk}.  We consider this
 problem from the viewpoint of the representation theory, being mainly interested in the \emph{representation type} (finite, tame or wild)
 and, for tame case, in the description of \emph{all} objects. So we do not deal with stable bundles and related topics, though something
 can be done in this direction too (cf. Section \ref{sec4}). We mostly consider algebras and varieties over an algebraically closed field $\Mk$,
 though some results remain valid in a more general setting.

  Recall a history of these investigations. The first general result concerning \emph{curve singularities} was obtained by H.\,Jacobinski
 and independently by A.\,Roiter and the author \cite{jac,dr}, who gave a criterion for a curve singularity to be of \emph{finite}
 Cohen--Macaulay type. We must note that there were no such \emph{words} as `curve singularity' or `Cohen--Macaulay module'
 in these papers; they were replaced by `commutative local ring of Krull dimension 1' and `torsion free module.' So the first paper,
 where these results were indeed related to curve singularities, was that of G.-M.\,Greuel and H.\,Kn\"orrer \cite{gk},
 where this criterion obtained the following wonderful form (see also \cite{yo}, which is a perfect survey on the Cohen--Macaulay finite type):

\smallskip
\emph{A curve singularity is \CM finite \iff it dominates a simple plane curve singularity} (in the sense of deformation theory, cf. 
 \cite{avg}). 

\smallskip
  This result was extended to \emph{hypersurface singularities} in \cite{kn,bgs}. Namely,

\smallskip
 \emph{A hypersurface singularity is \CM finite \iff it is simple (0-modal), i.e. of type $\rA_n,\,\rD_n\ (n\in\mN)$ or
 $\rE_n\ (n=6,7,8)$ } (cf. \cite{avg}). 

\smallskip
 At the same time H.\,\'Esnault and independently M.\,Auslander \cite{en,aus} (see also \cite{yo}) proved that

\smallskip
\emph{A surface singularity over an algebraically closed field $\Mk$ of characteristic 0 is \CM finite \iff it is a quotient singularity, i.e. is
 isomorphic to the ring of invariants $\Mk[[x,y]]^G$, where $G$ is a finite subgroup of $\mathbf{GL}(2,\Mk)$.}

\smallskip\rap
 By the way, we do not know whether this result has ever been generalized to the case of positive characteristic. 

  The first step towards \emph{tame} case was made by Schappert \cite{sch}, who proved that a plane curve singularity
 has at most 1-parameter families of ideals \iff it is \emph{strictly unimodal} \cite{wall}, or, the same, uni- or bimodal in the sense
 of \cite{avg}. In \cite{dg2} this result was extended to all curve singularities (one only has to replace in Schappert's theorem
`\emph{it is}' by `\emph{it dominates}'). Nevertheless, most of these singularities happened to be Cohen--Macaulay wild. Indeed,
 in \cite{dg1} \hbox{G.-M.}\,Greuel and the author showed that

\smallskip
\emph{A curve singularity, which is not \CM finite, is \CM tame \iff it dominates one of the singularities of type $\rT_{pq}$.}

\smallskip\rap
 In characteristic 0 these singularities (all of them are plane) are given by the equations $x^p+y^q+\la x^2y^2=0$ ($1/p+1/q\le1/2,\ \la\ne0$;
 if $(pq)=(36)\text{ or }(44)$ some special values of $\la$ must also be excluded). In positive characteristic the easiest way to define
 these singularities is by their parameterisations. 

 The case of surface singularities was first studied by C.\,Kahn \cite{kahn}. He proved that the so-called \emph{simple elliptic} singularities
 are \CM tame and described \CM modules over these singularities. Moreover, he elaborated a very general procedure that relates \CM
 modules over a normal surface singularity with vector bundles over the exceptional curve of its resolution. As for simple elliptic singularities
 the latter is an elliptic curve, he only had afterward to apply the classification of vector bundles over elliptic curves by M.\,Atiyah
 \cite{at}. By the way, till recently only vector bundles over a projective line \cite{gro} and over elliptic curves have been classified.

 To apply Kahn's technique to other surface singularities one has to know the corresponding results for vector bundles over projective curves.
 This investigation has been accomplished in \cite{dg3} with the following output:

\smallskip
 \emph{A projective curve is:
\begin{itemize}
\item 
  \VB\ finite if it is a configuration of projective lines of type $\rA_n$ (just $\mP^1$ if $n=1)$;
\item
  \VB\ tame if it is either an elliptic curve or a configuration of projective lines of type $\ti\rA_n$ (if $n=1$ it is a plane nodal cubic);
\item
  \VB\ wild otherwise.
\end{itemize}}

\smallskip\rap
 In finite and tame cases a complete description of vector bundles was obtained. These results (with the corresponding definitions and outline
 of proofs) are presented in the first part of the paper (Sections \ref{sec1}--\ref{sec8}). The technical background here is that of
 ``\emph{matrix problems},'' widely used before in analogous (and lots of other) questions. 

 The second part (Sections \ref{sec9}--\ref{sec14}) consists of applications to surface and hypersurface singularities (cf. \cite{dgk}).
 We recall the Kahn's reduction since it seems not well known to the audience. Moreover, we extend it to \emph{families} of \CM
 modules and vector bundles, which is necessary to deal accurately with tameness and wildness.
 Then we apply it, together with the results of the preceding sections, to the so-called \emph{minimally elliptic} singularities \cite{lau}.
 In this case a complete answer can be obtained:

\smallskip
 \emph{A minimally elliptic singularity is:
\begin{itemize}
\item 
 \CM tame if it is either simple elliptic or a cusp singularity;
\item
  \CM wild otherwise.
\end{itemize}}

 \smallskip\rap
 Here a \emph{cusp} singularity is a such one that  the exceptional curve of its minimal resolution is a configuration of projective lines
 of type $\ti\rA_n$ (the original definition by F.\,Hirzebruch \cite{hir} was different, though equivalent for the case $\Mk=\mC$). 
 For cusp singularities we get a complete description of \CM modules. By the way, it also gives possibility to fill up a flaw in the result
 for curve singularities, namely to give an explicit classification of \CM modules over the singularities of type $\rT_{pq}$. (Recall that in \cite{dg1}
 there was no such classification; their tameness was proved indirectly, using considerations from the deformation theory). We also obtain
 a description of \CM modules over the so-called \emph{log-canonical} singularities, since they are just quotients of simple elliptic or cusp
 by finite groups of automorphisms \cite{kaw}. There is some strong evidence that these cases are indeed the only tame ones, all other
 surface singularities being \CM wild.

 At last, we consider the case of \emph{hypersurface singularities}. Combining our results with the Kn\"orrer periodicity theorem
 \cite{kn,yo}, we obtain a classification of \CM modules over singularities \emph{of type} $\rT_{pqr}$, i.e. given by equations
 $x^p+y^q+z^r+\la xyz+Q(t_1,\dots,t_m)=0$, where $Q$ is a non-degenerate quadratic form, $1/p+1/q+1/r\le1,\ \la\ne0$
 (if $(p,q,r)=(2,3,6),(2,4,4)\text{ or }(3,3,3)$, some extra values of $\la$ must be excluded). 
 Again there is an evidence that all other hypersurface singularities are \CM wild.

 \section{An easy example: vector bundles on $\mP^1$}
 \label{sec1}

 A projective line $\mP^1$ is a union of two affine lines $\mA^1_i\ (i=0,1)$: if $(x_0:x_1)$ are homogeneous coordinates
 in $\mP^1,$ then $\mA^1_i=\setsuch{(x_0:x_1)}{x_i\ne0}$. The affine coordinate on $\mA^1_0$ is $t=x_1/x_0$
 and on $\mA^1_1$ $t^{-1}=x_0/x_1$. Thus we can identify $\mA^1_0$ with $\spec\Mk[t]$ and $\mA^1_1$
 with $\spec\Mk[t^{-1}]$; their intersection  is then $\spec\Mk[t,t^{-1}]$. Certainly, any projective module over
 $\Mk[t]$ is free, i.e. all vector bundles over an affine line are trivial. Therefore to define a vector bundle over $\mP^1$ one only has 
 to prescribe its rank $r$ and a \emph{gluing matrix} $A\in\GL(r,\Mk[t,t^{-1}])$. Changing bases in free modules over
 $\Mk[t]$ and $\Mk[t^{-1}]$ corresponds to the transformations $A\mps TAS$, where $S$ and $T$ are invertible matrices 
 of the same size, respectively over $\Mk[t]$ and over $\Mk[t^{-1}]$. Now an easy calculation, quite similar to that used in
 description of finitely generated modules over an euclidean ring, leads to the following
 \begin{lemma}\label{11}
  For any matrix $A\in\GL(r,\Mk[t,t^{-1}])$ there are matrices $S\in\GL(r,\Mk[t])$ and $T\in\GL(r,\Mk[t^{-1}])$ such
 that $SAT$ is a diagonal matrix\linebreak $\diag(t^{d_1},\dots,t^{d_r})$.
 \end{lemma} 
 Since $1\xx1$ matrix $(t^d)$ defines the line bundle $\kO_{\mP^1}(d)$, we get
 \begin{theorem}\label{12}
  Every vector bundle over a projective line uniquely decomposes into a direct sum of line bundles $\kO_{\mP^1}(d)$.

\em (As usually, `unique' in this context means that if $\kF\iso\bop_{i=1}^n\kL_i\iso\bop_{j=1}^m\kL'_j$, then $m=n$
 and there is a permutation $\tau$ of indices such that $\kL_i\iso\kL'_{\tau i}$ for all $i$.)
 \end{theorem} 

 This is a typical example of \emph{finite \VB\ type}: every indecomposable is a twist of one of them, namely $\kO_{\mP^1}$.
 I think that it is an important distinction between \emph{finite} and \emph{discrete} type. For instance, a quiver of type $\rA_\8$ 
 is a typical example of \emph{discrete}, but \emph{not finite} type: it has finitely many representation of each prescribed
 vector-dimension, but the dimensions of indecomposables can be arbitrary big.

 \section{A simple example: projective configurations of type $\rA$}
 \label{sec2}

 A \emph{projective configuration} is, by definition, a (singular, reduced) curve $X$ such that 
 \begin{enumerate}\brn
 \item 
  each irreducible component $X_i\ (i=1,\dots,s)$ of $X$ is \emph{rational}, i.e. its normalization is isomorphic to a
 projective line $\mP^1$;
 \item
  each singular point of $X$ is an \emph{ordinary double point}, i.e. a transversal intersection of two components (the latter may coincide,
 so it may be a self-intersection of a component).
 \end{enumerate}
 In particular, no three components pass through one point. For every projective configuration $X$ we define its \emph{intersection graph} (or \emph{dual graph}) $\De(X)$ as follows.
 \begin{itemize}
 \item 
  The \emph{vertices} of $\De(X)$ are irreducible components of $X$, or rather their indices $1,2,\dots,s$.
 \item
  The \emph{edges} of $\De(X)$ are singular points of $X$.
 \item
  An edge $p$ is incident to a vertex $i$ if $p\in X_i$; especially if $p$ is a self-intersection point of a component $X_i$,
 it gives rise to a \emph{loop} in the graph $\De(X)$.
 \end{itemize}
 This graph is non-oriented, but may contain loops and multiple edges between two vertices. Note that the graph $\De(X)$ does not, in general,
 define a projective configuration $X$ up to isomorphism. For instance, if $\De(X)$ is a graph of type $\ti D_3$, i.e. 
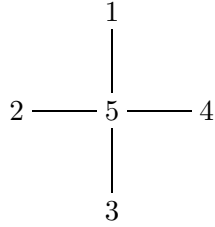
 $$
   \xymatrix{ & {1} \ar@{-}[d] & \\ {2}\ar@{-}[r] & {5}\ar@{-}[r] & {4} \\
  & {3}\ar@{-}[u]  }
 $$ 
 the position of 4 intersection points on the projective line $X_5$ corresponding to the central point depends on one parameter
 $\la\in\Mk\= \set{0,1}$: their harmonic ratio. For a fixed $\la$ these points can be chosen as $0,1,\la,\8$.

 We consider now the simplest case, when $\De(X)$ is of type $\rA_s$, i.e. 
 $$
   \xymatrix{
  {1}\ar@{-}[r] &{2}\ar@{-}[r]&{3}\ar@{-}[r]&{\cdots}\ar@{-}[r]&{s}
	} 
 $$ 
 Denote by $p_i\ (1\le i<s)$ the intersection point of $X_i$ and $X_{i+1}$. The calculation below is typical for the case of
 singular curves, though is the simplest example of this sort.

 The normalization $\tX$ of $X$ is just a disjoint union $\bsc_{i=1}^sX_i$. Each point $p_i$ gives rise to two points on
 $\tX$: $p'_i\in X_i$ and $p_i''\in X_{i+1}$. We may suppose that the isomorphisms $X_i\iso\mP^1$ are so chosen that
 $p'_i=\8$ and $p_i''=0$ (in homogeneous coordinates, respectively, $(0:1)$ and $(1:0)$). As the normalization mapping
 $\pi:\tX\to X$ is finite and birational, it induces an embedding $\kO=\kO_X\to\tO=\pi_*\kO_{\tX}$, and we can (and shall)
 identify any vector bundle $\kF$ over $\tX$ with its direct image $\pi_*\kF$. We denote by $\kJ$ the
 \emph{conductor} of $\tO$ in $\kO$, i.e. the biggest sheaf of $\tO$-ideals contained in $\kO$. In our example its sections
 are just those sections of $\tO$, which have zeros at all points $p_i'$ and $p_i''$ for each $i$. Note that
 $\kO/\kJ=\bop_{i=1}^{s-1}\Mk(p_i)$ and $\tO/\kJ\iso\bop_{i=1}^{s-1}(\Mk(p'_i)\xx\Mk(p''_i))$. 

 Let $\kG$ be a vector bundle over $X$ of rank $r$, $\tG=\tO\*_\kO\kG$. Then $\tG\sp\kG\sp\kJ\kG=\kJ\tG$.
 We already know that $\tG\iso\bop_{i=1}^s(\bop_{j=1}^r\kO_i(d_{ij}))$ for some integers $d_{ij}$, where $\kO_i=
 \kO_{X_i}$. For every $d$ 
$$
 \kO_i(d)/\kJ\kO_i(d)\iso\kO_i/\kJ\kO_i\iso
 \begin{cases}
  \Mk(p''_{i-1})\+\Mk(p'_i), &\text{if }\, 1<i<s, \\
  \Mk(p'_1),	&\text{if }\,i=1, \\
  \Mk(p''_{s-1}), &\text{if }\,i=s.
 \end{cases} 
 $$
 Therefore the factor $\tG/\kJ\tG$ is isomorphic to $\bop_{i=1}^{s-1}r(\Mk(p'_i)\+\Mk(p''_i))$. The factor
 $\kG/\kJ\kG$ is isomorphic to $\bop_{i=1}^sr\Mk(p_i)$, where each $r\Mk(p_i)$ is embedded into $r(\Mk(p'_i)\+
 \Mk(p''_i))$. Moreover, the projections of $r\Mk(p_i)$ onto both $r\Mk(p'_i)$ and $r\Mk(p''_i)$ are isomorphisms.
 On the contrary, given $r$-dimensional subspaces $V_i\in r(\Mk(p'_i)\+\Mk(p''_i))$ for each $i$ such that their projections
 onto both $r\Mk(p'_i)$ and $r\Mk(p''_i)$ are isomorphisms, we can construct a \VB\ $\kG$ over $X$ taking the preimage
 of $\bop_{i=1}^{s-1}V_i$ in $\tG/\kJ\tG$. Hence $\kG$ can be defined by a set of invertible $r\xx r$ matrices
 $\setsuch{M'_i,M''_i}{i=1,\dots,s-1}$ describing the projections of $V_i$ respectively onto $r\Mk(p'_i)$ and
 $r\Mk(p''_i)$. It is important that every row of these matrices has a \emph{weight} $d_{ij}$: the degree of the
 corresponding \VB\ $\kO_i(d_{ij})$ (this weight is common to the $j$-th rows of $M'_i$ and of $M''_{i-1}$).

 Certainly, we can change these matrices using automorphisms of $\tG$ and of $V_i$. Recall that 
$$
 \hom_{\kO_i}(\kO_i(d),\kO_i(d'))=
 \begin{cases}
   0, &\text{if }\, d>d',\\
   \Mk, &\text{if }\,d=d',\\
   \poly(d'-d), &\text{if }\,d<d',
 \end{cases} 
 $$
 where $\poly(m)$ is the space of homogeneous polynomials of degree $m$ over $\Mk$. Namely, if $s$ is a section of $\kO_i(d)$
 and $f\in\poly(d'-d)$, then $(fs)(\xi_0:\xi_1)=f(\xi_0,\xi_1)s(\xi_0:\xi_1)$. In particular, for any scalars $\la,\mu\in\Mk$ one
 can choose $f\in\poly(d'-d)$ such that $(fs)(0)=\la s(0)$ and $(fs)(\8)=\mu s(\8)$. Therefore two sets of matrices
 $\dM=\set{M'_i,M_i''}$ and $\dN=\set{N'_i,N''_i}$ define isomorphic \VB s over $X$ \iff $\dN$ can be obtained from $\dM$
 by a sequence of transformations of the following sorts:
 \begin{enumerate}\brn
  \item
   $M'_i\mps M'_iS$ and $M''_i\mps M''_iS$ for some $i$ and some invertible matrix $S$;
 \item
  $M'_i\mps T'M_i'$ and $M_{i-1}''\mps T''M_{i-1}''\ (1<i<s-1)$, where $T'=(t'_{jk})$ and $T''=(t''_{jk})$ are invertible
 matrices such that
       \begin{enumerate}
       \item   
         $t_{jk}'=t''_{jk}$ if $d_{ij}=d_{ik}$;
       \item
         $t_{jk}'=t''_{jk}=0$ if $d_{ij}<d_{ik}$;
       \end{enumerate}
 \item
  $M'_1\mps TM'_1$ for some invertible matrix $T$;
 \item
  $M''_{s-1}\mps TM''_{s-1}$ for some invertible matrix $T$.
  \end{enumerate}
 
 The following result is a rather simple exercise in matrix calculation.
 \begin{prop}\label{21}
  Using transformations \emph{(1--4)} from above one can transform any set $\dM=\set{M_i',M''_i}$ to the set 
 $\dI$ only consisting of unit matrices. 
 \end{prop} 
 (Actually, one has to start from $M_1'$, make it unit, then consider transformations of $M_1''$ that do not change this form of $M'_1$,
 etc.) 

 Evidently, it can be reformulated as a description of all \VB s over $X$.
 \begin{theorem}\label{22}
 Let $X$ be a projective configuration of type $\rA_s$.
\begin{enumerate}
\item 
  Every vector bundle over $X$ uniquely decomposes into a direct sum of line bundles.
 \item
  A line bundle $\kL$ over $X$ is uniquely determined by its \emph{vector-degree}, i.e. the sequence
 $\fD=\row ds$, where $d_i=\deg_{X_i}\kL$.
\end{enumerate}
 \end{theorem} 
 Especially every line bundle is a twist of the trivial line bundle $\kO$. Thus a projective configuration of type $\rA_s$, as well as
 projective line, is of  finite \VB\ type. Further we shall see that there are no more such curves. 

 Note that just the same calculation also gives a description of all \emph{torsion free} coherent sheaves $\kF$ over $X$.
 The distinctions are the following:
\begin{itemize}
\item 
  The ranks $r_j$ of the restrictions $\kF|X_j$ need not coincide, as well as the dimensions $m_i$ of the fibres $\kF(p_i)$.
 \item
  The projections $m_i\Mk(p_i)\to r_i\Mk(p'_i)$ and $m_i\Mk(p_i)\to r_{i+1}\Mk(p''_i)$ must be surjective (but not necessarily
 bijective).
 \item
  The mappings $m_i\Mk(p_i)\to r_i\Mk(p'_i)\+r_{i+1}\Mk(p''_i)$ must be injective.
\end{itemize}
 It means that the matrices $M'_i$ (of size $r_i\xx m_i$) and $M''_i$ (of size $r_{i+1}\xx m_i$) are not necessarily square,
 but have the ranks, respectively, $r_i$ and $r_{i+1}$, while the ``big'' matrix
 $$
   \mtr{M'_i\\M''_i}
 $$ 
 has rank $m_i$. In any case, using transformations (1--4) from above these matrices can be transformed to diagonal forms with
 1 and 0 on the diagonals. Thus an indecomposable torsion free sheaf is actually a vector bundle on a connected part of our
 configuration, i.e. on a curve consisting of components $X_k,X_{k+1},\dots,X_l$ for some $1\le k<l\le s$. Again
 we only get, up to twist, a finite number of indecomposables, so with this respect nothing changes.

 \section{Elliptic curves are \VB\ tame}
 \label{sec3}

 Another known case is that of \emph{elliptic curves}, i.e. smooth projective curves of genus 1. Such a curve can always be represented
 as a 2-fold covering of a projective line with 4 ramification points of degree 2, which can be chosen as $0,1,\la,\8\ \ (\la\in\Mk\=\set{0,1})$.
 If $\chr\Mk\ne2$ it can also be considered as a smooth cubic curve in $\mP^2$ given in one of its affine parts by the equation
 $y^2=x(x-1)(x-\la)$. Recall \cite[Section IV.4]{ha} that in this case the line bundles of a prescribed degree $d$ are in \oc\ with the
 points of the curve $X$. Namely, if we fix one point $o$, such a bundle is isomorphic to $\kO_X(x+(d-1)o)$ for a uniquely
 determined point $x$. Moreover, there is a line bundle $\kP$ on $X\xx X$ (\emph{Poincar\'e bundle}) such that, for every
 $x\in X$, 
 $$ 
  \kO_X(x+(d-1)o)\iso\kO_X(do)\*_{\kO_X}i_x^*\kP\iso i_x^*\kP(d(o\xx X)),
 $$
 where $i_x$ is the embedding $X\iso X\xx x\to X\xx X$. Thus the line bundles of degree $d$ form a 1-parameter family
 (parameterised by $X$).

 It so happens that the description of indecomposable vector bundles of an arbitrary rank and degree is quite similar. Nearby we present the
 results of Atiyah \cite{at} (with the modifications of Oda \cite{oda}, who has shown that the Atiyah's classification can be formulated in
 terms of families). We denote by $nx$ the closed subscheme of $X$ defined by the sheaf of ideals $\kO_X(-nx)$ and by
 $i_{nx}$ the embedding $X\xx nx\to X\xx X$.
 \begin{theorem}\label{31}
  For every pair of coprime integers $(r,d)$ with $r>0$ there is a \VB\ $\kP_{r,d}$ over $X\xx X$ such that every 
 indecomposable vector bundle over $X$ of rank $nr$ and degree $nd$, where $n$ is a positive integer, is isomorphic
 to ${p_1}_*i_{nx}^*\kP_{r,d}$ for a uniquely determined point $x\in X$. Moreover, $\kP_{r,d+mr}\iso
 \kP_{r,d}(m(o\xx X))$ and $\kP_{1,0}\iso\kP$.
 \end{theorem} 
  The proof of this theorem uses rather sophisticated considerations specific to elliptic curves, and we omit it referring to \cite{at,oda}.

 Theorem \ref{31} shows that any elliptic curve is \VB\ \emph{tame}: there are 1-paramenter families, at most one for any prescribed
 rank and degree, such that every indecomposable \VB\ over $X$ can be obtained from this family by specialization. The latter may
 include ``blowing,'' which means that we consider the values not only at points, but also at subschemes of the sort $nx$.
 (We shall give precise definitions in Section \ref{sec5}.)

 \section{Curves of genus $g>1$ are \VB\ wild}
 \label{sec4}

 Suppose now that $X$ is a smooth projective curve of genus $g>1$, $\kO=\kO_X$. Then for any two points
 $x\ne y$ from $X$ the \RR\ implies that $\hom_\kO(\kO(x),\kO(y))\iso\rH^0(X,\kO(y-x))=0$ and
 $\ext^1_\kO(\kO(x),\kO(y))\iso\rH^1(X,\kO(y-x))\ne0$. Fix 5 different points $x_1,\dots,x_5$ of the curve $X$, choose
 non-zero elements $\xi_{ij}\in\ext^1(\kO(x_j),\kO(x_i))$ for $i\ne j$ and consider vector bundles $\kF(A,B)$,
 where $A,B\in\mat(n\xx n,\Mk)$, and $\kF(A,B)$ is given as an extension
 $$
   0\larr \underbrace{ n(\kO(x_1)\+\kO(x_2))}_{\kB}\larr \kF(A,B)\larr
   \underbrace{n(\kO(x_3)\+\kO(x_4)\+\kO(x_5))}_{\kA} \larr 0
 $$ 
 corresponding to the element $\xi(A,B)$ of $\ext^1(\kA,\kB)$ presented by the matrix
 $$
   \mtr{ \xi_{13}I&\xi_{14}I&\xi_{15}I\\ \xi_{23}I&\xi_{24}A&\xi_{25}B  }
 $$ 
 ($I$ denotes the unit $n\xx n$ matrix). If $(A',B')$ is another pair of matrices, any homomorphism $\kF(A,B)\to\kF(A',B')$
 maps $\kO(x_i)$ to $\kO(x_i)$. It means that there are homomorphisms $\phi:\kA\to\kA'$ and $\psi:\kB\to\kB'$ such
 that $\psi\xi(A,B)=\xi(A',B')\phi$ (this is the Yoneda multiplication). Note that both $\phi$ and $\psi$ also map $\kO(x_i)$
 to $\kO(x_i)$ for each $i$. Now one can easily deduce that $\phi=\diag(S,S,S)$ and $\psi=\diag(S,S)$ for some
 matrix $S\in\mat(n'\xx n,\Mk)$ such that $SA=A'S$ and $SB=B'S$. 

 If we consider a pair $(A,B)$ as a representation of the free algebra $\Si_2$ in 2 generators, the correspondence $(A,B)\mps
 \kF(A,B)$ becomes a \emph{full, faithful, exact} functor $\md{\Si_2}\to\Vb(X)$. In particular, it maps non-isomorphic modules
 to non-isomorphic vector bundles and indecomposable modules to indecomposable vector bundles. Using terminology of the 
 representation theory of algebras, we say that the curve $X$ is \emph{\VB\ wild}. Again we give a precise definition in the next
 section.

 Recall that the algebra $\Si_2$ here can be replaced by \emph{any} finitely generated algebra $\La=\Mk\gnr{\lst am}$. Indeed,
 any $\La$-module $M$ such that $\dim_\Mk M=n$ is given by a set of matrices $\set{\lst Am}$ of size $n\xx n$. 
 One gets a full, faithful, exact functor $\md\La\to\md{\Si_2}$ mapping the module $M$ to the $\Si_2$-module of
 dimension $mn$ defined by the pair of matrices
 $$
   \mtr{\la_1I & 0 & \dots & 0\\ 0&\la_2I&\dots&0\\\hdotsfor4\\ 0&0&\dots &\la_nI},\quad
   \mtr{A_1&I&0&\dots &0\\ 0&A_2&I & \dots&0 \\\hdotsfor5\\ 0&0&0&\dots&A_n},
 $$ 
 where $\lst\la m$ are different elements from the field $\Mk$. Thus a classification of \VB s over $X$ would imply a
 classifications of \emph{all} representations of \emph{all} finitely generated algebras, the goal that perhaps nobody considers
 as achievable.
 
 As the first result of our investigation, we may formulate a theorem that describes \VB\ types of smooth projective curves.%
\footnote{The same results has been obtained by W.\,Scharlau (preprint of the M\"unster University). Moreover, he has also shown,
 almost in the same way, that every algebraic variety of dimension greater than 1 is \VB\ wild.}
 \begin{theorem}\label{41}
   A smooth projective curve $X$ is
\begin{itemize}
  \item 
   \VB\ finite if $X\iso\mP^1$;
\item
  \VB\ tame if it is an elliptic curve (i.e. of genus 1);
\item
  \VB\ wild otherwise.
\end{itemize}  
 \end{theorem} 

\begin{remk}\label{42}
  As $\kF(M)$ is an iterated extension of line bundles of degree 1, it is \emph{semi-stable} in the usual sense \cite{ses}. Thus even the
 classification of semi-stable vector bundles is wild in this case. The same can be shown in other cases too, though we do not bother to
 present explicit explanations.
\end{remk}

 \section{Vector bundle types: definitions}
 \label{sec5}

 The aim of this section is to precise the definitions concerning \VB\ types, especially make them available for non-smooth and even
 reducible curves. First we fix some notations.

 Let $X$ be a projective curve (connected, reduced, but maybe reducible) over an algebraically closed field $\Mk$. 
 We denote by $\lst Xs$ its irreducible components, by $\pi:\tX\to X$ its normalization, by $\sing X$
 the set of singular points, and by $\ting X=\pi^{-1}(\sing X)$ its preimage on $\tX$. Note that $\tX=\bsc_{i=1}^s\tX_i$,
 where $\tX_i$ is the normalization of $X_i$, so if $s>1$ it is not connected. We often write $\kO$ and $\tO$
 instead of, respectively, $\kO_X$ and $\kO_{\tX}$.

 Denote by $\Vb(X)$ the category of vector bundles over $X$.
 We always identify \emph{vector bundles} over $X$ with their sheaves of sections, thus with
 \emph{locally free coherent sheaves} of $\kO$-modules. If $\kF$ is such a sheaf, we set  
 $\tF=\kF\*_{\kO}\tO$; it is a vector bundle over $\tX$. As $\pi$ is finite and birational,  the direct image
 functor $\pi_*$ is a full embedding on the category of vector bundles, so we usually identify a vector bundle $\kG$ over
 $\tX$ with $\pi_*\kG$, which is a coherent sheaf on $X$ (but not a vector bundle over $X$). In particular, we usually
 identify $\tO$ with $\pi_*\tO$. 

 As $X$ is connected, every vector bundle $\kF$ has a constant rank $\rk\kF=\dim_\Mk\kF(x)$, where $x$ is an arbitrary
 closed point of $X$ and $\kF(x)=\kF_x/\gM_x\kF_x$. On the other hand, if $s>1$ we must consider the \emph{degree} of a vector
 bundle as a \emph{vector} $\deg\kF=\row ds$, where $d_i$ denotes the degree of the restriction $\kF|X_i$. The degree defines an
 epimorphism $\deg:\pic X\to\mZ^s$. We denote by $\pic^\circ X$ its kernel. We fix a section $s:\mZ^s\to\pic X$
 of this epimorphism and denote by $\kO(\fD)$ the line bundle $s(\fD) \ (\fD\in\mZ^s)$. Setting $\kF(\fD)=\kF\*_\kO\kO(\fD)$
 we define an action of the group $\mZ^s$ on vector bundles. 

 To define \VB\ tame and wild curves we need not individual sheaves, but their families, moreover, those with non-commutative bases.
 We provide the necessary definitions. Note that symbols like $\*,\ \hom,$ etc. always denote $\*_\Mk,\ \hom_\Mk,$ etc.
\begin{defin}\label{51}
  Let $X$ be a projective curve, $\La$ be a $\Mk$-algebra.
 \begin{enumerate}
\item 
 A \emph{family of \VB s over $X$ based on} $\La$ is a flat coherent sheaf of $\kO\*\La\op$-modules
 $\kF$ on $X$ (it is convenient to suppose that $\La$ acts on the right). We denote the category of such sheaves by $\Vb(X,\La)$.
\item
  Given such a family and any finite dimensional (over $\Mk$) $\La$-module $M$, we set $\kF(M)=\kF\*_\La M$;
  it is a vector bundle over $X$; moreover, for each vector $\fD\in\mZ^s$ we set $\kF(\fD,M)=\kF(\fD)\*_\La M$.
\\
 If $\La$ is commutative and $M=\Mk(x)=\La/\gM$, where $x$ is the closed point of $S=\spec\La$ corresponding to
 a maximal ideal $\gM\sb\La$, then $\kF$ can be considered as a family of vector bundles with the base $S$, and
 $\kF(M)=\kF(x)$ is the fibre of this family at the point $x$. If $S$ is connected (i.e. $\La$ is indecomposable), the 
 rank $\rk\kF(x)$ and the degree $\deg\kF(x)$ are constant on $S$; we call them the \emph{rank} and the \emph{degree
 of the family} $\kF$. If $M$ is an indecomposable, but not simple $\La$-module, $\kF(M)$ can be considered
 as a ``generalized'' fibre. For instance,
 if $M=\La/I$ for some ideal $I$, we consider $\kF(M)$ as the value of $\kF$ on the closed subscheme of $S$
 defined by the ideal $I$, just as we have done in Section \ref{sec3}. Note that we can consider families over arbitrary
 schemes, not only affine. The corresponding obvious changes in the definitions are left to the reader (cf. also \cite{dg3}).
\item
  A family $\kF$ of vector bundles over an algebra $\La$ is called \emph{strict} if, for every finite dimensional
 $\La$-modules $M,M'$,
  \begin{enumerate}
  \item  $\kF(M)\iso\kF(M')$ \iff $M\iso M'$;
  \item   $\kF(M)$ is indecomposable \iff so is $M$.
   \end{enumerate}
\item
  We call a curve $X$
  \begin{itemize}
  \item   
   \emph{\VB\ finite} if it has finitely many non-isomorphic indecomposable vector bundles \emph{up to twist}, i.e. there is a finite set
 of vector bundles $\lst\kF n$ such that every indecomposable vector bundle over $X$ is isomorphic to $\kF_k(\fD)$ for some
 $k\in\set{1,\dots,n}$ and some $\fD\in\mZ^s$;
  \item
   \emph{\VB\ tame} if there is a set $\kS$ of families of vector bundles over $X$
	satisfying the following conditions:
	\begin{enumerate}
	\item  every $\kF\in\kS$ is a strict family over a smooth connected curve $S_\kF$;
	\item
 the set  $\kS(r,\fD)=\setsuch{\kF\in\kS}{\rk\kF=r,\ \deg\kF=\fD}$  is finite  for each $r$ and $\fD$;
  	\item  for each $r$ and $\fD$ almost all (i.e. all but a finite number) indecomposable vector bundles over $X$ of
 	rank $r$ and degree $\fD$ are isomorphic to $\kF(\fC,M)$ for some $\kF\in\kS$, some vector $\fC\in\mZ^s$,
	 and some sky-scraper sheaf $M$, or, the same, a finite dimensional $\kO_{{S_\kF},x}$-module for some point $x\in S_\kF$
	\end{enumerate}
   (such a set $\kS$ is called a \emph{parameterising set} for vector bundles over $X$);
  \item
  \emph{\VB\ wild} if it possesses a strict family over every finitely generated $\Mk$-algebra $\La$.
  \end{itemize}
\item
  For a \VB\ tame curve $X$ and a parameterising set $\kS$ denote by $\nu(r,\fD,\kS)$ the cardinality of $\kS(r,\fD)$,
  and by $\nu(r,\fD,X)$ the smallest value of $\nu(r,\fD,\kS)$, when $\kS$ runs through all parameterising families. The curve $X$ is called
   \begin{itemize}
   \item  \emph{tame bounded} if there is a polynomial $N(r,\fD)$ such that $\nu(r,\fD,X)$ $\le N(r,\fD)$ for all $r$ and $\fD$;
   \item \emph{tame unbounded} otherwise.
   \end{itemize}
  For instance, elliptic curves are \VB\ tame bounded (actually in this case $\nu(r,d,X)=1$ for all $r$ and $d$).
\end{enumerate}
\end{defin}
\begin{remk}\label{52}
\begin{enumerate}
\item 
   The use of twists in the definitions of finite and tame is indeed indispensable. On the other hand, the parameterising families that we shall construct
 later actually cover \emph{all} vector bundles (up to twist), so the words `almost all' could be replaced by `all' in this context. Nevertheless, we
 have included them in order that our definition fits the usual pattern. 
\item
  One can be interested not only in vector bundles, but also in \emph{torsion free} sheaves, or even in all coherent sheaves.
 Certainly, all definitions of `finite,' `tame' or `wild' can be almost literally reproduced for these cases too. Moreover, the same calculations as
 for vector bundles show that nothing will change if we consider torsion free sheaves (though this time `almost all' is indispensable).
 We shall comment their structure at the corresponding places. The things become more complicated if we are interested in \emph{all}
 coherent sheaves over \emph{singular} curves, because they do not split into direct sums of torsion and torsion free ones.
 Nevertheless, I.\,Burban and the author have shown that even the \emph{derived category}
 of coherent sheaves remains tame for all \VB\ tame curves (cf. the talk of I.\,Burban presented at this workshop).
\item
  It is well-known that to prove that a curve $X$ is \VB\ wild it is enough to construct a strict family from $\Vb(X,\La)$, where $\La$
 is either the free algebra $\Mk\gnr{z_1,z_2}$ in two generators, or the polynomial algebra $\Mk[z_1,z_2]$, or the power series
 algebra $\Mk[[z_1,z_2]]$. Moreover, we can even replace them by a finite dimensional algebra, for instance,
 $\Mk[z_1,z_2]/(z_1^2,z_2^3,z_1z_2^2)$ or the path algebra of a wild quiver without cycles (the latter is even \emph{hereditary}
 that is sometimes convenient). In what follows we constantly use this observation.
\end{enumerate}
\end{remk}

 \section{Vector bundles and matrix problems}
 \label{sec6}

 We have already established the vector bundle types of smooth curves.
 Suppose now that $X$ is a singular curve. We use a procedure similar to that of Section \ref{sec2}.
 Namely, let $\kJ$ be the \emph{conductor} of $\tO$ in $\kO$,
 i.e. $\ann_\kO(\tO/\kO)$. Then $\kJ\tF=\kJ\kF\sb\kF\sb\tF$ for every vector bundle $\kF$. Denote by
 $\bF=\kO/\kJ$ and by $\tbF=\tO/\kJ$. Both $\bF$ and $\tbF$ have 0-dimensional support (respectively
 $\sing X$ and $\ting X$). Thus we may (and shall) identify them with the finite dimensional algebras of their sections
 $\Ga(X,\bF)$ and $\Ga(\tX,\tbF)$. Evidently $\bF$ is a subalgebra in $\tbF$. 

 We define the category $\fsm(X,\La)$, where $\La$ is a $\Mk$-algebra, as follows:
 \begin{itemize}
\item 
  Its \emph{objects} are the pairs $(\kA,M)$, where $\kA$ is an object from $\Vb(\tX,\La)$, i.e. a
 coherent flat sheaf of $\tO\*\La$-modules, and $M$ is a projective $\bF\*\La$-submodule in 
 $\kA/\kJ\kA$ such that the natural homomorphism $\tbF\*_\bF M\to\kA/\kJ\kA$ is an isomorphism.
\item
  A \emph{morphism} from $(\kA,M)$ to $(\kA',M')$ is a homomorphism $\phi:\kA\to\kA'$ such that 
 $\bar\phi(M)\sbe M'$, where $\bar\phi$ is the induced homomorphism $\kA/\kJ\kA\to\kA'/\kJ\kA'$.
\end{itemize}
 We write $\fsm(X)$ instead of $\fsm(X,\Mk)$. If $L$ is a finite dimensional $\La$-module and
 $\kP=(\kA,M)\in\fsm(X,\La)$, set $\kP(L)=(\kA\*_\La L,M\*_\La L)\in\fsm(X)$. Note also that the group $\mZ^s$ 
 acts on $\fsm(X)$. Indeed, the factors $\kA/\kJ\kA$ and $\kA(\fD)/\kJ\kA(\fD)$ are naturally isomorphic,
 so we can just set $\kP(\fD)=(\kA(\fD),M)$. It allows to transfer the definitions of finite, tame and wild types to the categories $\fsm(X)$
 (we leave it to the reader).

 Consider the functor $\bP:\Vb(X,\La)\to\fsm(X,\La)$ that maps a sheaf  $\kF$ to the pair $(\tF,\kF/\kJ\kF)$. Moreover, for
 each pair $\kP=(\kA,M)\in\fsm(X,\La)$ denote by $\bV(\kP)$ the preimage of $M$ in $\kA$; it is a sheaf of $\kO\*\La$-modules. 
 Obviously, any morphism $\phi:\kP\to\kP'$ induces a homomorphism $\bV(\kP)\to\bV(\kP')$ and one easily verifies the main property
 of these constructions:
 
\begin{prop}\label{p21}
  For every pair $\kP\in\fsm(X,\La)$ the sheaf $\bV(\kP)$ belongs to $\Vb(X,\La)$, and 
  the functors $\bP,\bV$ define an equivalence of the categories $\Vb(X,\La)$ and $\fsm(X,\La)$. Moreover, 
 $\bP(\kP(L))\iso\bP(\kP)(L)$ and $\bV(\kF(L))\iso\bV(\kF)(L)$ for every finite dimensional $\La$-module $L$.
\end{prop}
\begin{corol}\label{c22}
  A curve $X$ is \VB\ finite (tame, wild) \iff so is the category $\fsm(X)$.
\end{corol}

 It is easier to deal with the category $\fsm(X)$, because it can be identified with a bimodule category, so it carries us to the better explored
 world of ``matrix problems.'' (We refer to \cite{gr} and \cite{dg3} for the corresponding definitions and note that in \cite[Section 3]{dg3}
 the \emph{shifted bimodules} have been introduced, which are necessary for the application to vector bundles.)
  Namely, let $\sA=\Vb(\tX),\  \sB=\pro\bF$, the category of finitely generated projective right $\bF$-modules.
 Consider the $\sA\mbox{-}\sB$-bimodule $\sU$
 such that $\sU(P,\kA)=\hom_\bF(P,\kA/\kJ\kA)$, where $\kA\in\sA,\ P\in\sB$. In the category $\El(\sU)$ of elements of the
 bimodule $\sU$ (or matrices over $\sU$) consider the full subcategory $\El_c(\sU)$ consisting of all homomorphisms $\al:P\to\kA$
 such that the induced homomorphism $\tbF\*_\bF P\to\kA$ is an isomorphism. We call elements from $\El_c$ \emph{correct}.
\begin{prop}\label{p23}
  The categories $\fsm(X)$ and $\El_c(\sU)$ are equivalent.
\end{prop}
 Moreover, one can obviously define twists by $\fD\in\mZ^s$ on $\El(\sU)$ so that this equivalence is compatible with twists.
\begin{corol}\label{c24}
  A curve $X$ is \VB\ finite (tame, wild) \iff so is the category $\El_c(\sU)$.
\end{corol}

 \section{Vector bundle types: results}
 \label{sec7}

 Now we use the reduction to matrix problems to find the vector bundle types of singular curves. First we establish several \emph{wild} cases.
 We keep the notations and definitions of the preceding section and suppose that $X$ is indeed \emph{singular}. Since all proofs consist
 in more or less standard matrix calculations, we only write down the matrices that describe strict families, leaving the verification of
 strictness (always straightforward, though sometimes rather cumbersome) to the reader, who can also look into \cite{dg3}.

\begin{step}\label{s31}
  If one of the components $\tX_i$ is not rational, the curve $X$ is \VB\ wild. 
\end{step}
\begin{proof}
 Suppose that $\tX_1$ is of genus $g\ge1$. As $X$ is connected, there is a singular point $p$ that belongs to $X_1$.
 We suppose that $p$ has at least 2 preimages on $\tX$. If it only has one, the algebra $\tbF_p$ is not semi-simple, which
 simplifies the calculations. Let $\set{\lst pt}$ be all preimages of $p$, with $p_1\in \tX_1$, and let $Y$ be the component that contains
 $p_2$ (we allow $Y=\tX_1$). Let $\set{p_{t+1},\dots,p_l}$ be all other points from $\ting X$.
 Choose 4 different regular points $x_i\ (i=1,\dots,4)$ on $\tX_1$ and another regular point $y$ on $ Y$,
 and consider the family $\fU$ of elements of $\El_c(\sU)$ over the free algebra $\La=\Mk\gnr{z_1,z_2}$ that belongs
 to $\sU(\kA\*\La,4\bF\*\La)$, where $\kA=\bop_{k=1}^4 \tO(x_k+ky)$, all components of  $\fU$ in  
 $\hom_\bF(4\bF,\kA_{p_i}/\kJ\kA_{p_i})\*\La\ (1<i\le l)$ are unit matrices, and its component in
 $\hom_\bF(4\bF,\kA_{p_1}/\kJ\kA_{p_1})\*\La$ equals
 $$
   \mtr{0&0&0&1\\ 0&0&1&1\\ 0&1&1&z_1\\ 1&0&1&z_2}.
 $$ 
 Since $\hom_{\kO_{X_1}}(\kO_{X_1}(x_k),\kO_{X_1}(x_j))=0$ if $k\ne j$ and $\hom_{\kO_Y}(\kO_Y(ky),\kO_Y(jy))$
 $=0$ if $k>j$, one can check that $\fU$ is actually a strict family of correct elements. Hence $\El_c$ is wild and thus $X$ is
 \VB\ wild.
\end{proof}

 From now on we suppose that $\tX_i\iso\mP^1$ for every $i=1,\dots,s$.
\begin{step}\label{32}
  If the algebra $\tbF$ is not semi-simple, the curve $X$ is wild.
\end{step}
\rap (Note that $\tbF$ is semi-simple \iff all singular points of $X$ are \emph{ordinary multiple points}, i.e. such that at  each of them
 the number of linearly independent tangent directions to $X$ equals the multiplicity of this point.)
 \begin{proof}
 Choose a point $p\in\ting X$ such that $\tbF_p$ is not semi-simple, and a non-zero element $\al\in\tbF_p$ with $\ann_{\tbF_p}\al
 =\rad\tbF_p$. Let $\La$ be  the path algebra of the graph
 $$
 \begin{CD}
  5 @>z_5>> 4 @>z_4>> 3 \\
  && @Vz_3VV  @VVz_2V \\
  && 2 @>>z_1> 1
 \end{CD} 
 $$ 
 It is known to be wild \cite{df,na1}, so we only have to construct a strict family of correct elements over $\La$.
 Denote by $P_j$ the indecomposable
 projective $\La$-module corresponding to the vertex $j$ of the graph. We identify an arrow $z:j\to k$ with the corresponding 
 homomorphism $P_j\to P_k$ (left multiplication by $z$). Let $p\in\tX_1$. Choose a regular point $x$ on $\tX_1$. Set 
 \begin{align*} 
  \kA&=\tO\*(P_3\+P_1)\+\tO(x)\*(P_1\+P_2)\+\tO(2x)\*(P_3\+P_4)\+\tO(3x)\*P_5,\\
  B&=\bF\*(P_3\+P_1\+P_1\+P_2\+P_3\+P_4\+P_5) 
 \end{align*}  
 and consider the family $\fU$ of elements from $\El_c$ over $\La$ such that all components of $\fU$ except that in 
 $\hom_\bF(B,\kA_p/\kJ\kA_p)$ equal unit matrices, while the last is 
 $$ 
   \mtr{ &1&&0&0&0&\al\*1&0&0 \\  &0&&1&\al\*1&0&0&0&0 \\  &0&&0&1&\al\*z_1&\al\*z_2&0&0 \\ 
	&0&&0&0&1&0&\al\*z_3&0\\  &0&&0&0&0&1&\al\*z_4&0 \\ &0&&0&0&0&0&1&\al\*z_5 \\  
	&0&&0&0&0&0&0&1  } . 
 $$  
 Again a straightforward calculation shows that this family is strict, hence $X$ is \VB\ wild.
 \end{proof} 

\begin{step}\label{s33}
  If the curve $X$ has a singular point of multiplicity $m\ge3$, it is \VB\ wild.
\end{step}
 \begin{proof}
   Let $p$ be a point of multiplicity $l\ge 3$, which we suppose an ordinary multiple point, $\lst pl$ be its preimages on $\tX$.
 Denote by $Y_i$ the component of $\tX$ containing $p_i$ (some of them may coincide). Choose regular points $y_i\in Y_i$
 and set $\kA=\bop_{k=1}^4\tO(ky_1+ky_2+ky_3)$. Consider the family of elements from $\El_c(\sU)$ over
 $\La=\Mk\gnr{z_1,z_2}$ given by the element $\fU\in\hom_\bF(4\bF,\kA/\kJ\kA)$ that has unit matrices as all its
 components except those at the points $p_1$ and $p_2$, the last two being respectively
 $$
   \mtr{0&0&0&1\\ 0&0&1&0\\ 0&1&0&0\\ 1&0&0&0}	\ \text{ and }\ 
   \mtr{1&1&z_1&z_2 \\ 0&1&1&1 \\ 0&0&1&0 \\ 0&0&0&1}.
 $$ 
 One can check that $\fU$ is a strict family, so $X$ is \VB\ wild.
 \end{proof} 

 Thus from now on we only consider the case when all singular points of $X$ are \emph{ordinary double points} (or \emph{nodes}).
 It means that $X$ is a \emph{projective configuration} in the sense of Section \ref{sec2}, so its \emph{dual graph} $\De(X)$
 is defined.

\begin{step}\label{s34}
  If a vertex of $\De(X)$ is incident to three edges or to a loop and another edge, the curve $X$ is \VB\ wild.
\end{step}
\begin{proof}
 We consider the case when the graph $\De(X)$ is
 $$
 
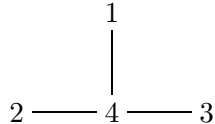
\xymatrix{
  & {1} \ar@{-}[d] \\
 {2}\ar@{-}[r] & {4}\ar@{-}[r] & {3}
 }
 $$
  (In other cases the calculations are even easier.) It means that the component $X_4$ intersects transversally the
 components $X_i\ (i=1,2,3)$ at the points $p_i\ (i=1,2,3)$. We denote by $p_{ij}$ the preimage of $p_i$ on
 the component $X_j\ (j=i\text{ or }j=4)$. If $\fU\in\hom_\bF(P,\kA/\kJ\kA)$ is an element of $\sU$, we denote by
 $\fU_{ij}$ the component of $\fU$ in $\hom_\bF(P,\kA_{p_{ij}}/\kJ\kA_{p_{ij}})$. We can choose homogeneous
 coordinates on $X_4\iso\mP^1$ so that $p_1=(1:0),\,p_2=(0:1),\,p_3=(1:1)$. Fix regular points 
 $x_i\in X_i$ and consider a family $\fU\in\hom_\bF({14}\bF,\kA)\*\Mk\gnr{z_1,z_2}$ of elements of $\sU$ such that
 \begin{align*}
       \kA&=\bop_{j=1}^{14} \tO(k_j(x_1+x_2+x_4)+l_jx_3)    \\
 \intertext{where} 
        k_j &=\begin{cases}
                  1 &\text{if }\ 1\le j\le3,\\ 
                  2 &\text{if }\ 4\le j\le7,\\
                  3 &\text{if }\ 8\le j\le11,\\
                  4 &\text{if }\ 12\le j\le14,
                \end{cases}           
 \end{align*}
        $ l_j=[(j+1)/2],\ \fU_{ij}$ is a unit matrix if  $ j\ne4\, \text{ or }\, i=1$,
 $$
        \fU_{24}=\mtr{ 0&\dots&0&1\\ 0&\dots&1&0\\ \hdotsfor4 \\ 1&\dots&0&0},
 $$
 while $\fU_{34}=(t_{pq})\ (p,q=1,\dots,14)$, where $t_{13,11}=z_1,\ t_{14,11}=z_2$, $t_{pq}=1$ if $p=q$
 or $(p,q)$ is from the set
 \begin{multline*}
 \{\,(5,3),(6,1),(7,2),(9,4),(10,6),(10,9),(11,5),(11,7),(11,10),\\(12,8),(13,9),(13,12),(14,10),(14,12)\,\},
 \end{multline*}
 and $t_{pq}=0$ otherwise.  Again a straightforward though cumbersome calculation shows that it is a strict family,
 thus $X$ is \VB\ wild.
  \end{proof}

 \begin{erem}
 Actually under the given shape of $\fU_{ij}$ for $(ij)\ne(34)$ the matrix $\fU_{34}$ splits into blocks
 $\fV_{kl}\ (1\le k\le4,\,1\le l\le7)$ (corresponding to the values  $k_j=k,\ l_j=l$), all of them with 2 columns,
  the number of rows is 3 for $k=1,4$ and 4 for $k=2,3$. With respect to the transformations that do not change
 other matrices these blocks form a representation of the pair of posets $(N,L)$, where $L$ is a chain with
 7 elements and $N$ is
 $$
  
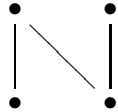
\xymatrix{
  {\bullet}\ar@{-}[d]\ar@{-}[dr] & {\bullet} \ar@{-}[d]\\
  {\bullet} & {\bullet}	}
 $$
  The matrix $\fU_{34}$ described above just presents a strict family of representations of this pair over
 $\Mk\gnr{z_1,z_2}$. Certainly, this matrix problem is known to be wild \cite{na2}, but we had to ensure
 the matrix $\fU_{34}$ to be invertible. That is why we had to take $L$ with 7 elements, though for the wildness of the 
 pair of posets 6 elements would suffice.  
 \end{erem} 

 Step \ref{s34} shows that  projective configurations that are not wild can only have the following dual graphs: 
  \begin{align*}
 &\text{chain}\ \rA_s& & 
\xymatrix{
  {1}\ar@{-}[r] &{2}\ar@{-}[r]&{3}\ar@{-}[r]&{\cdots}\ar@{-}[r]&{s}
	} \\ \intertext{or}
 &\text{cycle}\ \tilde\rA_s& &
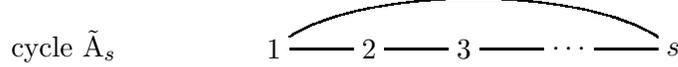
\xymatrix{
  {1}\ar@{-}[r] &{2}\ar@{-}[r]&{3}\ar@{-}[r]&{\cdots}\ar@{-}[r]&{s} \ar@{-}@(ul,ur)[llll]
	} 
  \end{align*}
 Here $s$ denotes the number of vertices. If $s=1$, there are no edges in the configuration of type $\rA_1$; thus the
 corresponding curve is just a projective line. The configuration $\tilde\rA_1$ corresponds to a rational curve with one ordinary
 double point (a nodal plane cubic, an affine part of which can be given by the equation $y^2=x^3+x^2$).

 We already know that the projective configurations of type $\rA_s$ are \VB\ finite.
 In the next section we show that all projective configurations of type $\tilde\rA_s$ are \VB\ tame
 unbounded, thus accomplishing the proof of the following theorem announced in the Introduction.

\begin{theorem}\label{t36}
  A projective curve is 
\begin{itemize}
\item  \VB\ \emph{finite} \iff it is a  projective configuration of type $\rA_s$,
\item \VB\ \emph{tame bounded} \iff it is a elliptic curve,
\item \VB\ \emph{tame unbounded} \iff it is a projective configuration of type $\tilde\rA_s$,
\item \VB\ \emph{wild} otherwise.
\end{itemize}
\end{theorem}
\rap
 Moreover, we present there a description of indecomposable vector bundles over the projective configurations of types $\tilde\rA_s$. 

 \section{Vector bundles over projective configurations of type $\ti\rA$}
 \label{sec8}

 Now we consider the projective configurations of type $\tilde\rA_s$. We follow the way of Section \ref{sec2} with evident changes.
 There are $s$ irreducible components $\lst\tX s$ of $\tX$ and $s$ singular points $\lst ps$ on $X$,
 each of them having two preimages $p_i'$ and $p_i''$ on $\tX$. We can arrange the numeration and coordinates on
 $\tX_i\iso\mP^1$ so that $p_i'=(1:0)\in \tX_i,\,p''_i=(0:1)\in\tX_{i+1}$ (we use the cyclic numeration modulo $s$, so
 $X_{s+1}=X_1$, etc.). Then $\bF=\prod_{i=1}^s\Mk_i$ and $\tbF=\prod_{i=1}^s(\Mk_i'\xx\Mk_i'')$, where
 $\Mk_i=\Mk(p_i),\,\Mk_i'=\Mk(p_i'),\,\Mk_i''=\Mk(p_i'')$. All these fields coincide with $\Mk$ and the embedding
 $\bF\to\tbF$ maps each $\Mk_i$ into $\Mk_i'\xx\Mk_i''$ diagonally. 

 Let $\fU\in\El_c(P,\kA)$. Then 
 $$
  P=r\bF,\quad \kA=\bop_{i=1}^s\bop_{k=1}^r\kO_i(d_{ik}),
 $$ 
 where $\kO_i=\kO_{\tX_i}$ and $d_{ik}\in\mZ$ are the degrees of direct summands; $r$ is the rank of $\kA$ (it must be
 constant). Note that $\kO_i(d)/\kJ\kO_i(d)\iso\Mk'_i\+\Mk''_{i-1}$ (again $\Mk_0=\Mk_s$). Choosing bases in
 each summand $\Mk'_i$ and $\Mk''_i$ as well as in each component $r\Mk_i$ of $M$, we present $\fU$
 as a set of $r\xx r$ invertible matrices $\setsuch{M_i',M_i''}{i=1,\dots,s}$. Moreover, the rows of these matrices are endowed
 with \emph{weights} $d_{ik}$ (it is the common weight of the $k$-th row of $M_i'$ and of $M_{i-1}''$). Taking into account
 the description of homomorphisms of vector bundles over $\mP^1$ from Section \ref{sec2}, we see that the automorphisms
 of $\kA$ and $P$ give rise to the following transformations of these matrices:
 \begin{enumerate}\brn
  \item
   $M'_i\mps M'_iS$ and $M''_i\mps M''_iS$ for some $i$ and some invertible matrix $S$;
 \item
  $M'_i\mps T'M_i'$ and $M_{i-1}''\mps T''M_{i-1}''\ (1\le i\le s)$, where $T'=(t'_{jk})$ and $T''=(t''_{jk})$ are invertible
 matrices such that
       \begin{enumerate}
       \item   
         $t_{jk}'=t''_{jk}$ if $d_{ij}=d_{ik}$;
       \item
         $t_{jk}'=t''_{jk}=0$ if $d_{ij}<d_{ik}$.
       \end{enumerate}
\end{enumerate}
 Two sets of matrices correspond to isomorphic vector bundles \iff one of them can be converted to the other by a sequence of
 such transformations.

 Note that this problem is much more complicated than that of Section \ref{sec2}. Even if $s=1$ and all rows have the same weight,
 it becomes the \emph{Kronecker problem}, or that of \emph{pencils of matrices}, which is of infinite type. Fortunately,
 one can recognize the arising matrix problem as belonging to the so-called ``representations of \emph{bunches of chains},'' or
 ``Gelfand problems,'' or ``clans'' (cf. \cite{cb,bon} or \cite[Appendix B]{dg3}; we refer to the last paper, because the presentation there is 
 more convenient to our purpose). Namely, in our case we have the pairs of chains 
 \begin{align*}
  \dE'_i&=\setsuch{ E'_{id}}{i=1,\dots,s;\ d\in\mZ},\quad \dF'_i=\set{f_i'},\\
  \dE''_i&=\setsuch{E''_{id}}{i=1,\dots,s;\ d\in\mZ},\quad \dF''_i=\set{f_i''}
 \end{align*} 
 with the natural ordering in each $\dE_i',\,\dE_i''$ (according to the index $d$),
 while the equivalence relation $\sim$ is given by the rules: 
 $$
  E'_{id}\sim E''_{i-1,d},\quad f'_i\sim f''_i.
 $$
 Slightly rearranging the list of indecomposable objects, we obtain the following result.

 \begin{theorem}\label{t41}
  Indecomposable vector bundles over a projective configuration of type $\tilde\rA_s$ are described by ``band data,''
 which are triples $(\fD,m,\la)$, where $m\in\mN,\ \la\in\Mk\=\set0$ and $\fD$ is a sequence from $\mZ^{rs}$ for some
 $r$ that is \emph{non-$s$-periodic}, i.e. cannot be presented as a repetition of a shorter sequence $\fC\in\mZ^{ls}\ (l<r)$.
 Two such triples $(\fD,m,\la)$ and $(\fD',m',\la')$ correspond to isomorphic vector bundles \iff $m=m',\,\la=\la'$ and
 $\fD'$ can be obtained from $\fD$ by an $s$-\emph{shift}, i.e. $\fD=\row d{rs}$, while
 $\fD'=(d_{ls+1},d_{ls+2}\dots,d_{rs},d_1,\dots,d_{ls})$ for some $l\le r$.
 \end{theorem} 
 
 Note that neither ``string data'' gives a set of \emph{invertible} matrices, though most of them can be interpreted as corresponding
 to \emph{torsion free}, but not locally free sheaves (cf. \cite{dg3} or Remark \ref{r43} below).

 Moreover, from the explicit description of indecomposable representations of a bunch of chains one can deduce an explicit description
 of vector bundles over such a configuration. Namely, the vector bundle $\kV=\kV(\fD,m,\la)$, where $\fD=\row d{rs}$,
 is a subsheaf of $\tV=\bop_{j=1}^{rs}m\kO_j(d_j)$ defined as follows: 
\begin{itemize} \label{vbd}
\item   
 If $x\notin\sing X$, $\kV_x=\tV_x$.
 \item
 Choose bases $\fE'_{jk}$ and $\fE''_{j,k}\ (j=1,\dots,rs,\,k=1,\dots,m)$ in each vector space $m(\kO_j(d_j)/\kJ\kO_j(d_j))_{p'_j}$
 and $m(\kO_{j+1}(d_{j+1})/\kJ\kO_{j+1}(d_{j+1}))_{p''_j}$ respectively, and set 
 $$
   \fE_{jk}=
 \begin{cases}
    \fE_{jk}'+\fE_{jk}'' &\text{if }\ 1\le j<rs,\\
    \fE_{rs,1}'+\la\fE_{rs,1}'' &\text{if }\ j=rs,\,k=1,\\
    \fE_{rs,k}'+\la\fE''_{rs,k}+\fE''_{rs,k-1} &\text{if }\ j=rs,\,k>1.
 \end{cases} 
 $$ 
 \item
  For each singular point $p_i$ the stalk $\kV_{p_i}$ is generated by the preimages of $\fE_{jk}$ with $j\equiv i\pmod s$.
\end{itemize}
 Especially $\rk\kV(\fD,m,\la)=mr$ and $\deg\kV(\fD,m,\la)=\row\de s$, where $\de_i=\sum_{j\equiv i\pmod s}d_j$.
 For instance, if $s=1$ (the case of nodal cubic) and $m=1$, one can present this vector bundle as the following gluing of line bundles
 over $\tX=\mP^1$:  
 $$ 
   
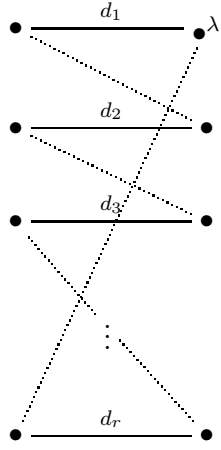
\xymatrix{
  {\bul} \ar@{-}[rr]^{d_1}\ar@{.}[rrd] && {\bul^\la}\ar@{.}[ddddll] \\
  {\bul} \ar@{-}[rr]^{d_2}\ar@{.}[rrd] && {\bul} \\
  {\bul} \ar@{-}[rr]^{d_3}\ar@{.}[rd] && {\bul} \\
   & \vdots\ar@{.}[rd] & \\
  {\bul} \ar@{-}[rr]^{d_r} && {\bul}
	} 
 $$  
  Here horizontal lines symbolize line bundles over $\tX$ of the superscripted degrees, their left (right) ends are basic elements of these bundles
 at the point $0=(1:0)$ (respectively $\8=(0:1)$), and the dotted lines show which of them must be glued. All gluings are trivial, except that going from
 the uppermost right point to the lowermost left one, where we glue one vector to another multiplied by $\la$. If $m>1$, one has to take $m$
 copies of each \VB\ from this picture, make again trivial all gluings except the last one, where identifications must be made using the Jordan
 $m\xx m$ cell with eigenvalue $\la$. The necessary changes for $s>1$ are quite obvious.

 \begin{corol}\label{c42}
  Projective configurations of type $\tilde\rA_s$ are \VB\ tame unbounded.
 \end{corol} 
 \begin{proof}
  Let $\La=\Mk[x,x^{-1}]$.  For each vector $\fD\in\mZ^{rs}$ we define a family $\kV_\fD$
 of vector bundles of rank $r$ with base $\spec\La=\mA^1\=\set0$ just as we have defined the vector bundles $\kV(\fD,1,\la)$,
 but replacing $\kO_j$ by $\kO_j\*\Mk[x,x^{-1}]$ and $\la$ by $x$. Obviously, then
 $\kV(\fD,m,\la)\iso\kV_\fD(L)$, where $L=\La/(x-\la)^m$, hence $X$ is \VB\ tame. As we allow twists, it is enough to
 consider sequences $\fD$ of non-negative integers (even those containing 0). Evidently, for every fixed degree $\de
 =\row\de s$ and every fixed rank $r$ there are finitely many such $\fD$ with $\sum_{j\equiv i\pmod s}d_j=\de_i$.
 On the other hand, the number of such sequences $\fD$ grows exponentially when $r\to\8$. Thus $X$ is indeed tame unbounded.
 \end{proof} 

 \begin{remk}\label{r43}
  Just in the same way we can describe \emph{torsion free} sheaves over any configuration of type $\tilde\rA_s$. The difference is that it is allowed
 not to glue all basic vectors at $0$ and $\8$, leaving the first and the last ones ``free'' (for instance, in the picture above we can exclude
 one of the gluings). Note that in this case we can make all gluings trivial.
 It means that there are no families of torsion free sheaves that are not vector bundles; all of them stand apart. See \cite{dg3} for details.
 \end{remk} 

 \section{\CM modules: generalities}
 \label{sec9}

 From now on we consider \CM modules over \emph{normal surface singularities}. For our purpose, it means a complete noetherian
 algebra $\bA$ over an algebraically closed field $\Mk$ such that
 \begin{itemize}
 \item 
    $\bA$ is \emph{normal}, i.e. integral and integrally closed;
 \item
   $\mathrm{Kr.dim}\,\bA=2$;
 \item
  $\bA/\gM\iso\Mk$, where $\gM$ denotes the maximal ideal of $\bA$.
 \end{itemize}
  Note that for 2-dimensional integral rings `\emph{Cohen--Macaulay}' means that $\bA=\bap_{\he\gP=1}\bA_\gP$, while `\emph{normal}'
 means that  $\bA_\gP$ is  a \dvr\ for each prime ideal $\gP$ of height 1. Moreover, a finitely generated module $M$ over such a ring
 is \CM \iff it is torsion free and $M=\bap_{\he\gP=1}M_\gP$. We denote $M\vv=\hom_\bA(M,\bA)$; it is always a \CM $\bA$-module.
 For \CM $\bA$-modules $M,N$ we denote by $M\bot_\bA N$ their ``reflexive product'' $(M\*_\bA N)\vvv$.

 A \emph{family of \CM $\bA$-modules} based on an algebra $\La$ is defined as a finitely generated
 $\bA\mbox{-}\La$-bimodule $\kM$ such that
\begin{itemize}
\item   $\kM$ is flat as $\La$-module;
 \item
  for every finite dimensional $\La$-module $L$ the $\bA$-module $\kM(L)=\kM\*_\La L$ is Cohen--Macaulay.
\end{itemize}
 Obviously, the latter condition is only to be checked for simple modules $L$. Having this notion, one can define
 \emph{strict families, \CM tame and wild singularities} just as it has been done in Section \ref{sec5} for vector bundles,
 so we omit the details of these definitions. We also leave to the reader the obvious changes in these definitions,
 when families based on arbitrary $\Mk$-schemes are considered. We call  modules $\kM(L)$ (generalized)
 fibres of $\kM$.

   We denote by $S=\spec\bA$ and by $p$ the unique closed point of $S$ (corresponding to the ideal $\gM$).
 Recall that there always is a \emph{resolution} of such a singularity, i.e. a projective birational morphism of schemes
 $\pi:X\to S$, where $X$ is regular, such that the restriction of $\pi$ onto $\oX=X\=\pi^{-1}(p)$ is an
 isomorphism $\oX\to\oS=S\=\set p$, cf. e.g. \cite{lip}. The reduced preimage $E=\pi^{-1}(p)_{\mathrm{red}}$ is
 called the \emph{exceptional curve} of this resolution. It is indeed a projective curve, possibly singular and even reducible.
 We denote by $\lst Es$ its irreducible components. We always identify $\oX$ and $\oS$ so that the diagram  
  $$ 
 \begin{CD}
  \oX @>i>> X\\
  @|	@VV\pi V \\
  \oS @>>j>S
 \end{CD} 
 $$  
 where $i,\,j$ are embeddings, commutes.

 To obtain a criterion of \CM finiteness, as well as for some other results, the following considerations are important.
 Suppose that a finite group $G$ acts on a normal surface singularity $\bB$ and $\bA=\bB^G$  is the ring of invariants.
 It is again a normal surface singularity. A \CM $\bB$-module $N$ is called \emph{induced from}
 $\bA$ if it is isomorphic to $\bB\bot_\bA M$ for some
 \CM $\bA$-module $M$. We say that $\bB$ is \emph{unramified in codimension} 1 if, for every prime ideal $\gP\sb\bA$ of height 1,
 $\bB_\gP/\gP\bB_\gP$ is a separable algebra over the field $\bA_\gP/\gP\bA_\gP$. Equivalently, the natural epimorphism of
 $\bB$-bimodules $\bB_\gP\*_{\bA_\gP}\bB_\gP\to\bB_\gP$ splits.  For instance it is so if $G$ acts freely on the set
 of prime ideals of $\bB$ of height 1.

 \begin{prop}[cf. \cite{her}]\label{p53} 
  If the order $g=|G|$ is invertible in $\Mk$, every \CM $\bA$-module is a direct summand of a \CM $\bB$-module considered as $\bA$-module.
 If, moreover, $\bB$ is unramified in codimension 1, every \CM $\bB$-module is a direct summand of a module induced from $\bA$.
 \end{prop} 
 \begin{proof}
   Let $M$ be any \CM $\bA$-module. We identify it with the $\bA$-submodule $1\*M$ of $\bB\*_\bA M$.
 There is a retraction of $\bA$-modules $\bB\*_\bA M\to M$ mapping $b\*v$ to $g^{-1}\sum_{\si\in G}\si(b)\*v$.
 Thus $\bB\*_\bA M\iso M\+M'$ for some $M'$. Taking second duals, we get $\bB\bot_\bA M\iso M\+(M')\vvv$.

 Suppose now that $\bB$ is unramified in codimension 1. Denote by $\bK$ and $\bL$ respectively the fields
 of fractions of $\bA$ and $\bB$. Then $\bL$ is a Galois extension of $\bK$ with Galois group $G$. Hence
 $\bL^e=\bL\*_\bK\bL=\bop_{\si\in G}\bL^\si$ as $\bL$-bimodule, where
 $\bL^\si=\setsuch{\la\in\bL}{\forall\, {\al\in\bL}\,\ \al\la=\la\si(\al)}$ 
 \cite{dk}. Especially there is a unique element $\eps\in\bL^e$ such that $\la\eps=\eps\la$ and $\phi(\eps)=1$, where $\phi$
 denotes the natural epimorphism $\bL^e\to\bL$. Since the restriction of $\phi$ onto $\bB^e_\gQ$ splits, $\eps\in\bB_\gQ^e$.
 Therefore $\eps\in\bap_{\he\gQ=1}\bB^e_\gQ=\bB\bot_\bA\bB$, so the natural epimorphism $\bB\bot_\bA\bB\to\bB$
 splits too. Let $N$ be any \CM $\bB$-module. Then $\bB\*_\bA N\iso(\bB\*_\bA\bB)\*_\bA N$, hence
 $\bB\bot_\bA N\iso(\bB\bot_\bA\bB)\bot_\bB N$, which gives a natural epimorphism
 $\bB\*_\bA N\to N$. It arises from $\bB\bot_\bA\bB\to\bB$, thus splits, and $N$ is a direct summand of $\bB\bot_\bA N$.
 \end{proof} 

 Obviously, in the situation of Proposition \ref{p53} an indecomposable \CM $\bA$-module (respectively $\bB$-module) is isomorphic to
 a direct summand of an indecomposable \CM $\bB$-module (respectively $\bA$-module).

 \begin{corol}\label{c54}
  If the order $|G|$ is invertible in $\Mk$ and $\bB$ is \CM finite, tame or wild, then so is $\bA$. If, moreover, $\bB$ is unramified
 in codimension 1 and $\bA$ is \CM finite, tame or wild, then so is $\bB$.
 \end{corol} 

 It implies immediately the \'Esnault--Auslander criterion. Namely call $\bA$ a \emph{quotient singularity} if it is isomorphic to a ring of invariants
 $\bR^G$, where $\bR=\Mk[[x,y]]$ and $G$ is a finite group of automorphisms of $\bR$ (it is well-known that in this case
 we may always suppose that $G$ acts linearly, so it is a finite subgroup of $\mathrm{GL}\,(2,\Mk)$).

 \begin{theorem}[\'Esnault--Auslander] \label{t54}
  Suppose that $\chr\Mk=0$. Then $\bA$ is \CM finite \iff it is a quotient singularity.
 \end{theorem} 
 \begin{proof}
  \emph{Necessity.}  If $\bA$ is \CM finite, it has in particular finitely many rank 1 \CM modules, or, the same, divisorial ideals (up to isomorphism).
 Thus its Picard group (group of classes of divisors) is finite. It is known (cf. e.g. \cite{lip}) that such a singularity is \emph{rational}, which means that
 $\rH^1(X,\kO_X)=0$. Consider the \emph{canonical divisor} $K=K_S$, i.e. the class of a (rational) differential 2-form. It is of finite order
 in the Picard group, i.e. $nK$ is a principal divisor for some $n$. Choose an ideal $J$ of class $K$; then $J^{-n}=\al\bA$ and one can
 consider the ring $\bB=\bop_{i=0}^{n-1}J^{-i}t^i$ with $t^n=\al^{-1}$. It is Gorenstein and $\bA\iso\bB^H$, where $H$
 is the cyclic group of order $n$ acting naturally on $\bB$, namely leaving elements of $\bA$ intact and mapping $t$ to $\eps t$,
 where $\eps$ is a primitive root of 1. For every prime $\gQ\sb\bA$ of height 1 $\bA_\gQ$ is a \dvr, hence $J_\gQ=\ga\bA$ for
 some $\ga$ such that $\al=\mu\ga^{-n}$ with an invertible $\mu$. Moreover, since $\chr\Mk=0$, $\mu$ is an $n$-th power of some element,
 so we may suppose that $\al=\ga^{-n}$ and $\bB_\gQ=\bA_\gQ[\ga ^{-1}t]\iso\bA_{\gQ}[x]/(x^n-1)\iso\bA_{\gQ}^n$. Therefore
 $\bB$ is unramified in codimension 1, so it is also \CM finite by Proposition \ref{p53}, thus rational. But all rational Gorenstein
 singularities are well-known \cite{lau}. They are \emph{rational double points}, or \emph{du Val singularities}. All of them are 
 quotient singularities. Therefore $\bA$ is quotient too.

 \emph{Sufficiency} follows directly from Proposition \ref{p53}, since any regular local ring is \CM finite: all \CM modules are free.
 \end{proof} 

 \begin{remk}\label{r55}
  As far as I know, the finiteness criterion is still unknown if $\chr\Mk>0$, though it seems very plausible that the answer must be the same
 (maybe modulo some minor changes of definitions).
 \end{remk} 
 
 \section{Kahn's reduction}
 \label{sec10}

 In this section we recall the main results of the Kahn's paper \cite{kahn} and extend them to families of \CM modules and \VB s.

An \emph{(exceptional) cycle} on $X$ is a divisor $Z=\sum_{i=1}^sz_iE_i$.
 If $z_i\ge0$ it is called \emph{effective}. We treat an effective cycle as a projective curve (non-reduced if $z_i>1$
 for some $i$), namely the subvariety of $X$ defined by the sheaf of ideals $\kO_X(-Z)$. We also denote by
 $\om_X$ the \emph{dualizing sheaf} of $X$ and by $\om_Z=\om_X(Z)\*_{\kO_X}\kO_Z$ the dualizing sheaf
 of $Z$. The latter defines the \emph{Serre's duality} 
 $$
   \rH^i(E,\kF)\iso\rD\rH^{1-i}(E,\kF\vv\*_{\kO_Z}\om_Z) \qquad (i=0,1)
 $$ 
 for every vector bundle $\kF$ on $Z$, where $\rD V=\hom(V,\Mk)$, the dual vector space to $V$, and
 $\kF\vv=\Hom_{\kO_Z}(\kF,\kO_Z)$.

 For any coherent sheaf $\kF$ on $X$ denote by $\kF^g$ the image of the evaluation mapping $\Ga(X,\kF)\*\kO_X\to\kF$.
 If $\kF^g=\kF$, we say that $\kF$ is \emph{generated by global sections}, or \emph{globally spanned}. If the support of the
 factor $\kF/\kF^g$ is 0-dimensional (i.e. a finite set of closed points), we say that $\kF$ is \emph{generically generated}
 by global sections, or \emph{generically spanned}.

 The main notion of the Kahn's theory is the following.
 \begin{defin}\label{d61}
  An effective cycle $Z$ is called a \emph{weak reduction cycle} if
\begin{enumerate}\brn
\item 
 the sheaf $\kO_Z(-Z)$ is generically spanned;
 \item
  $\rH^1(E,\kO_Z(-Z))=0$.\\
 It is called a \emph{reduction cycle} if, moreover,
 \item
  the sheaf $\om_Z\vv$ is generically spanned.
\end{enumerate}
 \end{defin} 
  A reduction cycle always exists: it easily follows from the fact that the \emph{intersection form} $(A.B)$ is negative definite on the
 group of all exceptional cycles \cite{lip}. 

 We identify $\bA$-modules with quasi-coherent sheaves over $S$ (their sheafifications). In particular, we consider the inverse
 image functor $\pi^*:\mdl\bA\to\coh X$. Even if $M\in\Cm(\bA)$, usually $\pi^*M$ can have torsion. So we define the 
 functor $\pi\vs:\Cm(\bA)\to\Vb(X)$ setting $\pi\vs(M)=(\pi^*M)\vvv$, where $\kF\vv=\Hom_{\kO_X}(\kF,\kO_X)$.
 As $\pi$ is isomorphism outside $E$ and any \CM module $M$ is completely defined by its stalks outside $p$,
 $M\iso\pi_*\pi\vs M$, so this functor is full and faithful. The following theorem describes its image.

 \begin{theorem}[\cite{kahn}]\label{t62}
  A \VB\ $\kF$ over $X$ is isomorphic to $\pi\vs M$ for some \CM $\bA$-module $M$ \iff
\begin{enumerate}\brn
\item
  $\kF$  is generically spanned;
 \item
  the restriction $\Ga(X,\kF)\to\Ga(\oX,\kF)$ is surjective.
\end{enumerate}
 \rm We call such vector bundles \emph{full} and denote by $\Vb^f(X)$ the subcategory of full vector bundles.
 \end{theorem} 

 If $Z$ is any effective cycle, we can consider the functor ``restriction on $Z$,'' $\res_Z:\Vb(X)\to\Vb(Z)$ mapping $\kF$
 to $\kF/\kF(-Z)$. 
 \begin{theorem}[\cite{kahn}]\label{p62}
  \begin{enumerate}
 \item
  The functor $\res_Z$ is \emph{dense}, i.e. every vector bundle over $Z$ is isomorphic to $\res_Z\kF$ for some vector bundle
 $\kF$ over $X$.
 \item
  If $Z$ is a reduction cycle, the restriction of $\res_Z$ onto $\Vb^f(X)$ 
 maps non-isomorphic vector bundles to non-isomorphic ones.
\\ {\rm We call a \VB\ over $Z$ \emph{full} if it is isomorphic to $\res_Z\kF$, where $\kF\in\Vb^f(X)$, and denote by
 $\Vb^f(Z)$ the category of full vector bundles over $Z$.}
 \item
  The functor  $R_Z=\res_Z\circ\pi^\sharp$ induces a representation equivalence\linebreak $\Cm(\bA)\to\Vb^f(Z)$.
\end{enumerate}
 \end{theorem} 
  Note that this functor cannot be faithful, since $\hom$-spaces in the category $\Vb(Z)$ are finite dimensional. Moreover, it can
 map indecomposable vector bundles to decomposable ones, cf. Theorem \ref{t71} below.

 Kahn also gives a description of full vector bundles over a weak reduction cycle.
 
 \begin{theorem}[\cite{kahn}]\label{p63}
  If $Z$ is a weak reduction cycle, the following conditions on a vector bundle $\kE\in\Vb(Z)$are equivalent:
 \begin{enumerate}\brn
 \item
  $\kE$ is full.
 \item
    $\kE$ is generically spanned and
  there is a vector bundle $\kE_2$ over the cycle $2Z$ such that $\kE_2|Z\iso\kE$ and the mapping
 $\rH^0(\kE(Z))\to\rH^1(\kE)$ induced by the canonical exact sequence $0\to\kE\to\kE_2(Z)\to\kE(Z)\to0$
 is injective.  \end{enumerate} 
\rm  The latter sequence is obtained by tensoring 
 with $\kE(Z)$ the exact sequence $0\to\kO_Z(-Z)\to\kO_{2Z}\to\kO_Z\to 0$.
 \end{theorem} 
 
 Actually we need a generalization of  Theorems~\ref{t62},\ref{p63} for \emph{families} of vector bundles and modules
 based on an algebra $\La$.

 \begin{defin}\label{d63}
   Let $\pi:X\to S$ be a resolution of a normal surface singularity and $\La$ be a $\Mk$-algebra (maybe non-commutative).
\begin{itemize}
 \item
  For any coherent sheaf of $\kO_X\*\La$-bimodules $\kF$
   we denote by $\kF^g$ the image of the evaluation mapping $\Ga(X,\kF)\*(\kO_X\*\La)\to\kF$ and say that
  $\kF$ is \emph{generically spanned} if the support of $\kF/\kF^g$ is 0-dimensional, i.e. a finite set of closed points.
\item 
  We call a family of \VB s $\kF$ over $X$ based on $\La$ \emph{full} if it is isomorphic to $\pi\vs\kM=(\pi^*\kM)\vvv$, where
 $\kM$ is a family of \CM modules over $S$ based on $\La$ and $\kN\vv=\Hom_{\kO_X\*\La\op}(\kN,\kO_X\*\La\op)$
 for any family $\kN$ of $\kO_X$-modules based on $\La$.
 \item
  We call a family of vector bundles over an effective cycle $C$ \emph{full} if it possess a full lifting to a family of \VB\ 
 over $X$.
\end{itemize}
 \end{defin} 

 To extend Kahn's results to families, one needs some restrictions on the base algebra $\La$. For our purpose it is enough to
 consider algebras $\La$ such that $\gdim\La\le2$. The advantage is that in this case any kernel of a mapping between
 two flat $\La$-modules is also flat. Especially, if $\kF$ is a family of $\kO_X$-modules based on
 $\La$, $U\sbe X$ is any open subset and $U=\bup_iU_i$ is its affine open covering, then $\Ga(U,\kF)$ is the kernel of 
 the natural mapping $\bop_i\Ga(U_i,\kF)\to\bop_{i,j}\Ga(U_i\cap U_j,\kF)$, thus $\La$-flat. As a corollary, if
 $\pi:X\to Y$ is a proper morphism, the direct image $\pi_*\kF$ is $\La$-flat, so is a family of $\kO_Y$-modules based
 on $\La$. (We need `proper' in order  $\pi_*\kF$ to be coherent.) The same is true if we consider families based on regular
 schemes of dimension at most 2, since their local rings are of global dimension at most 2.

 \begin{theorem}\label{t64}
    Suppose that $Z$ is a weak reduction cycle for a resolution $\pi:X\to S$ of a normal surface singularity and 
 $\gdim\La\le2$.
 \begin{enumerate}
 \item 
  Let $\kF$ be a family of \VB s over $X$ based on $\La$ such that
     \begin{enumerate}
     \item 
      $\kF$ is generically spanned;
     \item
       the restriction $\Ga(X,\kF)\to\Ga(\oX,\kF)$ is surjective.
     \end{enumerate}
  Then $\kF$ is full.
 \item
    Let $\kE$ be a family of vector bundles over $Z$ such that
   \begin{enumerate}
   \item   
     $\kE$ is generically spanned;
    \item
    there is a lifting of $\kE$ to a \VB\ over $2Z$ such that
     the induced mapping $\rH^0(E,\kE(Z))\to\rH^1(E,\kE)$ is injective.
   \end{enumerate}
  Then $\kE$ is full.
 \end{enumerate}
 \end{theorem} 
 \begin{proof}
   1. Set $\kM=\pi_*\kF$, $\kF'=\pi^*\kM/(\kO_X\text{-torsion})$. Then $\kF'$ can be considered as a subsheaf of $\kF$ containing
 $\kF^g$ (since global sections of $\kM$ are the same as those of $\kF$). Note that $\kM(L)\iso\pi_*\kF(L)$ is a \CM module by
 \cite[Proposition 6.3.1]{ega}. Therefore $\kM$ is a family of \CM modules. Since $\kF$ is generically generated, its stalks $\kF_x$
 coincide with $\kF^g_x$ provided $x$ is not a closed point. But as $X$ is normal, any family $\kN\vv$ is completely determined by
 its stalks at non-closed points. Thus $\kF\iso\kF\vvv\iso(\kF')\vvv\iso\pi\vs\kM$.

 2. Suppose that $\kE_n$ is a family of vector bundles over $nZ$ such that $\kE_n\*_{\kO_{nZ}}\kO_Z\iso\kE$.
 Then the obstruction for lifting $\kE_n$ to a family of vector bundles over $(n+1)Z$ lies in $\rH^2(E,\Hom(\kE,\kE(-n)))$.
 (It can be shown just in the same way as in \cite{pet}, where the complex analytic case was considered.) But this cohomology
 space is $0$ since $\dim E=1$. Hence such a lifting is always possible, so we can construct a sequence $\kE_n$, where each
 $\kE_n$ is a family of vector bundles over $nZ$, such that $\kE_n\*_{\kO_{nZ}}\kO_{(n-1)Z}\iso\kE_{n-1}$. Taking inverse limit
 we get a family $\kF$ of vector bundles over $X$ based on $\La$ such that $\kE_1=\kE$ and
 $\kF\*_{\kO_X}\kO_{nZ}\iso\kE_n$ for all $n$. 
 We only have to show that $\kF$ is full. According to the Theorem on Formal Functions \cite[Theorem III.11.1]{ha}, the $\gM$-adic
 completion of $\rH^1(X,\kF(-Z))$ coincides with $\varprojlim\rH^1(E,\kE_n(-Z))$. The natural exact sequence
 $$
   0\to\kO_{Z}(-nZ)\larr \kO_{(n+1)Z} \larr \kO_{nZ}\to 0  
 $$
 tensored with $\kF$ gives an exact sequence
 \begin{equation}\label{seq}
    0\to\kE(-nZ)\larr \kE_{n+1} \larr \kE_n\to 0.
 \end{equation} 
 Since both $\kE$ and $\kO_Z(-Z)$ are generically spanned, so is $\kE(-nZ)$ for all $n\ge0$. It means that for every $n>0$ there is a
 homomorphism $m\kO_Z\to\kE(-(n-1)Z)$ with the cokernel support of dimension 0. Twisting it by $-Z$ we get a homomorphism
 $m\kO_Z(-Z)\to\kE(-nZ)$ with the same property, which induces an epimorphism $m\,\rH^1(E,\kO_Z(-Z))\to\rH^1(E,\kE(-nZ))$.
 Since $Z$ is a weak reduction cycle, we get $\rH^1(E,\kE(-nZ))=0$ for all $n>0$. The exact sequence \eqref{seq} twisted by $-Z$
 gives an exact sequence of cohomologies
 $$
    \rH^1(E,\kE(-(n+1)Z))\larr \rH^1(E,\kE_{n+1}(-Z))\larr \rH^1(E,\kE_n(-Z))\larr 0,
 $$ 
 wherefrom we can deduce that $\rH^1(E,\kE_n(-Z))=0$ by an obvious induction. Therefore $\rH^1(X,\kF(-Z))=0$ and the exact
 sequence $0\to\kF(-Z)\to\kF\to\kE\to0$ shows that every global section of $\kE$ can be lifted to a global section of $\kF$. But
 outside $E$ every quasicoherent sheaf over $X$ is generated by its global sections. Hence $\kF$ is generically spanned, i.e. the condition
 1(a) holds. Note also that the equality $\rH^1(X,\kF(-Z))=0$ implies that $\rH^1(X,\kF)\iso\rH^1(E,\kE)$.
 
  To verify the condition 1(b) we use local cohomologies \cite{gro2}, especially the exact sequence 
 $$
      \rH^0(X,\kF)\larr\rH^0(\oX,\kF)\to\rH^1_E(X,\kF)\to\rH^1(X,\kF),
 $$ 
 which shows that the condition 1(b) can be reformulated as follows:

\smallskip
{}\quad 1(b$'$) \emph{the mapping $\rH^1_E(X,\kF)\to\rH^1(X,\kF)$ is injective}.

\smallskip\noindent
 Due to \cite[Lemma B.2]{wahl} we can identify $\rH^1_E(X,\kF)$ with $\varinjlim\rH^0(E,\kE_n(nZ))$, where
 the limit is taken along homomorphisms $$\mu_n:\rH^0(E,\kE_{nZ}(nZ))\to\rH^0(E,\kE_{n+1}((n+1)Z))$$
 arising from the natural exact sequence
 $$
    0\to\kO_{nZ}(-Z)\larr \kO_{(n+1)Z} \larr \kO_{Z}\to 0   
 $$ 
 tensored by $\,\kE_{n+1}((n+1)Z)\,$. Especially all $\mu_n$ are injective. We shall show that under our conditions they are
 also surjective, or, equivalently, all homomorphisms
 $
   \rH^0(E,\kE((n+1)Z))\larr \rH^1(E,\kE_n(nZ))
 $ 
 are injective. Actually we shall prove that even their compositions with the restrictions $\kE_n(nZ)\to\kE(nZ)$, i.e.
 homomorphisms
 $$
   \de_n:\rH^0(E,\kE((n+1)Z))\larr\rH^1(E,\kE(nZ))
 $$ 
 are injective. The latter arise from the exact sequence
 $$
   0\to\kE(nZ)\larr \kE_2((n+1)Z) \larr \kE((n+1)Z)\larr 0,     
 $$ 
 thus, by the condition 2(b), we may suppose that it is injective for $n=0$. Since $Z$ is a weak reduction cycle, the sheaf
 $\kO_X(-Z)$is generically spanned, hence all sheaves
 $\kO_X(-nZ)$ are generically spanned too, so there is a homomorphism $m\kO_X\to\kO_X(-nZ)$ with the cokernel supported
 on a finite set of closed points. Then the dual mapping $\kO_X(nZ)\to m\kO_X$ is a monomorphism. Tensoring with $\kE(Z)$,
 we get a monomorphism $\phi:\kE((n+1)Z)\to m\kE(Z)$, wherefrom we get the following commutative diagram: 
 $$ 
 \begin{CD} 
  \rH^0(E,\kE((n+1)Z) @>\de_n>> \rH^1(E,\kE(nZ)) \\
  @V\rH^0(\phi)VV	@VVV \\
  m\,\rH^0(E,\kE(Z)) @>{m\,\de_0}>> \rH^1(E,\kE).
 \end{CD} 
 $$ 
 Since both $\de_0$ and $\rH^0(\phi)$ are injective, so is $\de_n$. 

 Thus $\rH^1_E(X,\kF)\iso\rH^0(E,\kE(Z))$. Since also $\rH^1(X,\kF)\iso\rH^1(E,\kE)$, we see that
 the condition 1(b) is actually equivalent to the condition 2(b). It accomplishes the proof.
 \end{proof} 

 We also need the following important, though rather simple, observation.

 \begin{prop}\label{p65}
  If  a family $\kF$ of vector bundles over $X$ based on $\La$ is full, so are also all its fibres $\kF(L)=\kF\*_\La L$.
 \end{prop} 
 \begin{proof}
  First show that $\rH^1(\oX,\kF)=0$. Note that since $E$ is a closed subscheme of a regular scheme $X$, it can be locally
 defined by one equation  \cite[Proposition II.6.11]{ha}.
 Therefore $\rH^2_E(X,\kF)=0$ \cite{gro2} and the mapping $\rH^1(X,\kF)\to\rH^1(\oX,\kF)$ is surjective.
  The exact sequence $0\to\kF(-E)\to\kF\to\kE\to0$, where $\kE=\kF\*_{\kO_X}\kO_E$, together with the obvious
 equalities $\rH^i(\oX,\kE)=0$, implies that $\rH^i(\oX,\kF)\iso\rH^i(\oX,\kF(-nE))$ for all $n$. But $X$ is projective over
 the affine scheme $S$ and $-E$ is ample, so $\rH^1(X,\kF(-nE))=0$ for some $n$, hence also $\rH^1(\oX,\kF)=0$.
  
 Now from the K\"unneth formulae \cite{ce} we obtain a commutative diagram 
 $$ 
 \begin{CD}
  \rH^0(X,\kF)\*_\La L @>>> \rH^0(X,\kF\*_\La L)  \\
  @VVV @VVV\\
 \rH^0(\oX,\kF)\*_\La L @>>> \rH^0(\oX,\kF\*_\La L) ,
 \end{CD} 
 $$ 
 where the lower horizontal arrow is an isomorphism and the left vertical arrow is surjective. Hence the right vertical arrow is surjective
 too, which means that $\kF(L)$ is full. 
 \end{proof} 

 \begin{remk}\label{r66}
 Obviously, the subcategory $\Vb^f(X)\sb\Vb(X)$ is closed under direct summands. On the contrary, it is not the case for
 the subcategory $\Vb^f(Z)\sb\Vb(Z)$ (cf. for instance Theorem \ref{t71} below). 
 The same is true for families too. That is why, even under the conditions of Theorem \ref{t64}.2,
 we cannot claim that a full lifting $\kF$ of $\kE$ is strict if so is $\kE$, though it is true that $\kF(L)\iso\kF(L')$ implies 
 $L\iso L'$. On the other hand, if $\kF$ is strict, so is the family $\kM=\pi_*\kF$ of \CM $\bA$-modules, since the restriction
 of $\pi_*$ onto $\Vb^f(X)$ is an equivalence $\Vb^f(X)\str\sim\to\Cm(\bA)$.
 \end{remk}

 \section{\CM types: minimal elliptic case}
 \label{sec11}

  Recall that the \emph{fundamental cycle}  of a resolution $\pi:X\to S$ of a normal surface singularity is the smallest effective
 cycle $Z$ such that $(Z.E_i)\le0$ for each irreducible component $E_i$ of the exceptional curve $E$ \cite{lau}.
 This singularity is called \emph{minimally elliptic} \cite{lau,reid} if it is Gorenstein and
 $\nH^1(\kO_X)=1$ for some (hence for any) resolution $\pi:X\to S$. If $\pi$ is minimal, an equivalent condition
 is: $(Z+K.E_i)=0$ for all $i$, where $Z$ is the fundamental cycle of this resolution and $K$ is the canonical
 divisor of $X$. In particular, then $\om_Z\iso\kO_Z$. One can easily check that $Z$ is a reduction cycle in this case
 \cite{kahn,dgk}.

 For minimally elliptic singularities the criterion of Theorem \ref{p63}  can be essentially simplified and restated in terms of
 $\kE$ itself, without references to liftings. 

 \begin{theorem}[\cite{kahn}]\label{t71}
  Let $\kE$ be a \VB\ over $Z$. It is full \iff $\kE\iso\kG\+n\kO_Z$, where
  \begin{enumerate}\brn
\item   $\kG$ is generically spanned;
 \item
  $\rH^1(\kG)=0$;
 \item
  $n\ge\nH^0(\kG(Z))$.
\end{enumerate}
 Under these conditions the full lifting $\kF$ of $\kE$ is indecomposable \iff $\kG$ is indecomposable and either $\kE=\kO_Z$ or 
 $\kG\not\iso\kO_Z$ and $n=\nH^0(\kG(Z))$.
 \end{theorem} 

 We can extend this result to families as follows.

 \begin{theorem}\label{t72}\brn
  Let $\bA$ be a minimally elliptic singularity, $\pi:X\to S$ be its minimal resolution, and $Z$ be the fundamental
 cycle of this resolution (which is known to be a reduction cycle). Suppose that $\kG$ is a family
 of \VB s over $Z$ based on an algebra $\La$ with $\gdim\La\le2$ such that
  \begin{enumerate}
  \item 
  $\kG$ is generically spanned;
 \item
  $\rH^1(E,\kG)=0$;
\end{enumerate}
 We set $P_0=\rH^0(E,\kG(Z))$. Let also $P$ be a projective $\La$-module such that there is an embedding $P_0\to P$.
 Then the family $\kE=\kG\+(\kO_Z\* P)$ is full.
 \end{theorem} 
 \begin{proof}
   Obviously $\kE$ is generically spanned. Moreover, 
 \begin{align*}
  \rH^1(E,\kE)&=\rH^1(E,\kO_Z\*P)\iso \rH^1(E,\kO)\*P\iso P,\\
  \rH^0(E,\kE(Z))&=\rH^0(E,\kG(Z))=P_0,
 \end{align*} 
 since $\rH^0(E,\kO_Z(Z))=0$ for any exceptional cycle $Z$ (cf. e.g. \cite[Chapter 4, Exercise 14]{reid}). 
 We have already seen that there is a lifting $\kE'$ of $\kE$ to a family of vector bundles over $2Z$.
 It gives an exact sequence
 $$
   0\larr \kE(-Z)\larr \kE'\larr \kE\larr0
 $$ 
 Denote by $\xi$ the corresponding element of  
 $\ext^1_{\kO_{2Z}\*\La}(\kE,\kE(-Z))$ and by $\de$ the induced mapping $\rH^0(E,\kE(Z))\to\rH^1(E,\kE)$.
 One can easily see that any element from $\ext^1_{\kO_{2Z}\*\La}(\kE,\kE(-Z))$ that is of the form $\xi+\eta$,
 where
 $
 \eta\in\ext^1_{\kO_Z\*\La}(\kE,\kE(-Z))\iso\ext^1_{\kO_Z\*\La}(\kE(Z),\kE),
 $
 also defines a lifting of $\kE$ to a family of vector bundles.
 Moreover, such an element induces the mapping $\de+\de_\eta:\rH^0(E,\kE(Z))\to\rH^1(E,\kE)$, where
 $\de_\eta(s)=\eta s$, the Yoneda product of $\eta\in\ext^1_{\kO_Z\*\La}(\kE(Z),\kE)$ with 
 $s\in\rH^0(E,\kE(Z))\iso\hom_{\kO_Z\*\La}(\kO_Z\*\La,\kE(Z))$. Choose $\eta$ from
 \begin{multline*}
  \ext^1_{\kO_Z\*\La}(\kG(Z),\kO_Z\*P)\iso\ext^1_{\kO_Z\*\La}(\kO_Z\*D_\La P,\kG\vv(-Z))
 \iso\\ \iso\rH^1(E,\kG\vv(-Z)\*_\La P),
 \end{multline*}
 where $\kG\vv=\Hom_{\kO_Z\*\La}(\kG,\kO_Z\*\La),\ D_\La P=\hom_\La(P,\La)$. Due to the K\"unneth formulae, 
 $\rH^1(E,\kG\vv(-Z)\*_\La P)\iso\rH^1(E,\kG\vv(-Z))\*_\La P$. The Serre's duality implies that
 $\rH^1(E,\kG\vv(-Z))\iso D_\La\rH^0(E,\kG(Z))\iso D_\La P_0$.
 Since $D_\La P_0\*_\La P\iso\hom_\La(P_0,P)$, the mapping $(s,\eta)\to\eta s$ is just the evaluation homomorphism
 $P_0\xx\hom_\La(P_0,P)\to P$. Thus $\de_\eta$ coincides with $\eta$ as a homomorphism $P_0\to P$.
 Therefore we can choose $\eta$ so that $\de+\de_\eta$ becomes any prescribed homomorphism $P_0\to P$,
 for instance an embedding. Then the corresponding lifting $\kE_2$ of $\kE$ to a family of vector bundles over $2Z$ satisfies the
 condition (2b) of Theorem \ref{t64}, hence $\kE$ is full.
 \end{proof} 

 \begin{corol}\label{c73}\brn
  Suppose that the conditions of Theorem \ref{t72} hold, as well as the following:
 \begin{enumerate}
 \setcounter{enumi}2
 \item  
  $\kG$ is strict and has no fibres isomorphic to $\kO_Z$;
 \item
  $\rH^1(E,\kG(Z))$ is flat as $\La$-module;
 \item
  $P=P_0$.
\end{enumerate}
 Then the full lifting $\kF$ of the family $\kE$ to a family of vector bundles over $X$ is also strict, as well as the
 family $\kM=\pi_*\kF$ of \CM modules over $\bA$.
 \end{corol} 
 \begin{proof}
  Let $L$ be an indecomposable $\La$-module. The K\"unneth formulae imply that
 $$
  \rH^0(E,\kG(Z)\*_\La L)\iso\rH^0(E,\kG(Z))\*_\La L=P\*_\La L
 $$
  since $\rH^1(E,\kG(Z))$ is flat, and
 $$
  \rH^1(E,\kG\*_\La L)\iso\rH^1(E,\kG)\*_\La L=0
  $$ 
 since $\rH^2(E,\kG)=0$.
 Hence $\kG(L)$ satisfies the indecomposability conditions from Theorem \ref{t71} (with $n=\dim_{\Mk}P\*_\La L$).
 Since $\kF(L)$ is full and $\res_Z\kF(L)\iso \kE(L)$, $\kF(L)$ is indecomposable. Moreover, if $\kF(L)\iso\kF(L')$, then
 $\kG(L)\+n\kO_Z\iso\kG(L')\+n'\kO_Z$ for $n'=\dim_{\Mk}P\*_\La L'$. As neither $\kG(L)$ nor $\kG(L')$ has
 direct summands isomorphic to $\kO_Z$, the Krull--Schmidt theorem for vector bundles implies that $\kG(L)\iso\kG(L')$,
 thus $L\iso L'$, so $\kF$ and hence $\kM$ are strict.
 \end{proof} 
 
 Now we can use the results on vector bundle types to define \CM types of minimally elliptic singularities. Recall some definitions.
 \begin{defin}\label{d74}
  A normal surface singularity with a minimal resolutions $\pi:X\to S$ and exceptional curve $E$ is called 
 \begin{itemize}
\item 
  \emph{simple elliptic} if $E$ is a smooth elliptic curve \cite{sai};
 \item
  \emph{cusp} if $E$ is a projective configuration of type $\ti\rA$.
\end{itemize}
 The latter is not the original definition (in the case $\chr\Mk=0$), but is equivalent to it \cite{hir,kar}. We accept it as a \emph{definition}
 of cusp singularities in the case $\chr\Mk>0$ too.
 It is easy to see that both simple elliptic and cusp singularities are minimally elliptic; moreover the fundamental cycle in these cases coincides
 with the exceptional curve $E$.
 \end{defin} 

 Note that, according to the \'Esnault--Auslander Theorem \ref{t54}, neither minimally elliptic singularity can be \CM finite
 (it also follows directly from the next theorem, even if $\chr\Mk>0$).
 
 \begin{theorem}\label{t75}
  A minimally elliptic singularity $\bA$ is
 \begin{itemize}
\item 
 \CM \emph{tame bounded} if it simple elliptic;
 \item
  \CM \emph{tame unbounded} if it is a cusp singularity;
 \item
  \CM \emph{wild} otherwise.
\end{itemize}
 \end{theorem} 
 \begin{proof}
   Suppose first that $\bA$ is neither simple elliptic nor cusp. Then the exceptional curve $E$, hence also the fundamental cycle $Z$
 is \VB\ wild, i.e. possess a strict family $\kG$ of vector bundles over the polynomial ring $\bR=\Mk[x,y]$. Replacing $\kG$ by
 $\kG(m)$ for $m$ big enough we may suppose that $\kG$ is generically (even globally) spanned and $\rH^1(E,\kG)=0$.
 Set $P=\rH^0(E,\kG(Z))$ and $Q=\rH^1(E,\kG(Z))$. 
 There is an element $g\in\bR$ such that $Q[g^{-1}]$ is flat over $\bR[g^{-1}]$. Moreover, there is at most one point $z\in\mA^2$ 
 such that $\kG(\Mk(z))\iso\kO_Z$. Choose $x\in\mA^2$ such that $x\ne z$ and $g(x)\ne0$. Let $\gX$ be the corresponding
 maximal ideal of $\bR$. Set $\La=\hat\bR_{\gX}$ (the $\gX$-adic completion) and $\hat\kG=\kG\*_\bR\La$. Then $\hat\kG$
 is a family of vector bundles over $Z$ based on $\La$ that satisfied conditions of Corollary \ref{c73}. Thus the family
 $\hat\kG\+\kO_Z\*\hat P_\gX$ has a full lifting $\kF$ to $X$, which is strict, and the family of \CM $\bA$-modules
 $\kM=\pi_*\kF$ based on $\La$ is strict too. Since $\La\iso\Mk[[x,y]]$, $\bA$ is \CM wild (cf. Remark \ref{52}.3).

 If $\bA$ is simple elliptic,  we can use the Atiyah--Oda description of vector bundles over $E$ together with the calculations from
 \cite{at}. The latter give the following values of cohomologies for vector bundles $\kP_{r,d}(nx)={p_1}_*i^*_{nx}\kP_{r,d}$ from
 Theorem \ref{31}:
 \begin{align*}
  \nH^0(E,\kP_{r,d}(nx))&=
 \begin{cases}
   nd, &\text{ if } d>0,\\
   1, &\text{ if } d=0 \text{ and } x=o,\\
   0 &\text{ otherwise}.
 \end{cases} \\
 \nH^1(E,\kP_{r,d}(nx))&=
 \begin{cases}
   nd, &\text{ if } d<0,\\
  1, &\text{ if } d=0 \text{ and } x=o,\\
   0 &\text{ otherwise}.       
 \end{cases} 
 \end{align*} 
 In particular, $\kP_{r,d}(nx)$ is generically spanned \iff either $d\ge r$ or $d=0,\,r=n=1$, $x=o$ (the latter gives the trivial bundle $\kO_E$),
 and $Q=\rH^1(E,\kP_{r,d}(Z))$ is flat for $d\ne br$, where $b=-(E. E)$
 (recall that in this case $Z=E$, hence $\kP(Z)=\kP(-bo)$). If $d=br$ the sheaf $Q$ is no more flat, since its fibre at the 
 point $o$ jumps, but its restriction onto $E^o=E\={o}$ is flat. Denote by $\kP^o_{r,d}$ the restriction of $\kP_{r,d}$ onto $E^o$
 and define the families $\kE_{r,d}$ of vector bundles over $E$, where $r\le d,\,(r,d)=1$, as follows:
 $$
   \kE_{r,d}=
 \begin{cases}
  \kP_{r,d}, &\text{ if } r\le d< br,\\
  \kP_{r,d}\+\kO_E\*(d-br)\kO_E, &\text{ if } d>br,\\
  \kP^o_{1,b}, &\text{ if } r=1,\,d=b. 
 \end{cases}  
 $$ 
 In the first two cases this family is based on $E$, in the last one it is based on $E^o$. These families satisfy the conditions of Corollary \ref{c73},
 hence they can be lifted to families $\kF_{r,d}$ of vector bundles over $X$. Denote $\kM_{r,d}=\pi_*\kF_{r,d}$. Then Theorems \ref{31}
 and \ref{t72} directly imply a description of \CM $\bA$-modules (also obtained in \cite{kahn}). Below we denote by $a^+$ and
 $a^-$ respectively the \emph{positive} and \emph{negative part} of a number $a$ defined as
 $$
   a^+=
 \begin{cases}
  a &\text{ if } a>0,\\ 0 &\text{ otherwise};
 \end{cases} \qquad
 a^-= \begin{cases}
  -a &\text{ if } a<0,\\ 0 &\text{ otherwise}.
 \end{cases} 
 $$ 

 \begin{theorem}\label{t76}
  If $\bA$ is simply elliptic, all indecomposable \CM $\bA$-modules are:
 \begin{itemize}
 \item
  $\bA$, of rank $1$;
\item 
 $\kM_{r,d}(nx)$, of rank $n(r+(d-br)^+)$, where $r\le d,\, (r,d)=1$ and $x\in E$, $x\ne o$ if $\,r=1,\,d=b$;
 \item
 $M_n$, of rank $n+1$, where $R_E M_n\iso\kP_{1,b}(no)\+\kO_E$.
\end{itemize}
 In particular $\bA$ is \CM tame bounded, since the number of families $\kM_{r,d}$ with $r+(d-br)^+=m$ is 
 at most $b\vi(m)$, where $\vi(m)$ is the Euler function.
 \end{theorem} 

 Let now $\bA$ be a cusp singularity, $E=\bup_{i=1}^sE_i$ be the exceptional curve of its minimal resolution,
 $E_i$ being the irreducible components so arranged that $E_i\cap E_{i+1}\ne\0$ (as before, we set $E_{s+i}=E_i$),
 $b_i=-(E_i.E_i)$, $\fB=\row bs$ and $\fB^r$ be the $r$-fold repetition of $\fB$, i.e.
 $$
  \fB^r=(\lst bs,\lst bs,\dots,\lst bs)\quad (r\, \text{ times}).
  $$
 In this case $Z=E$ and $\kG(Z)\iso\kG(-\fB)$.
 If $\kG$ is a vector bundle over $E$, $\ti\kG$ is its lifting to the
 normalization $\ti E$, one can calculate cohomologies of $\kG$ using the long exact sequence of cohomologies corresponding to the
 exact sequence $0\to\kG\to\ti\kG\to\ti\kG/\kG\to0$ and the known values of $\rH^i(\ti E,\ti\kG)$ (recall that $\ti E$ is just a
 union of projective lines). Applying this procedure to the vector bundles $\kV(\fD,m,\la)$ from page \pageref{vbd}, one gets that these
 values can be calculated as follows (cf. \cite{dgk} for details). Define a \emph{positive part} of $\fD=\row d{rs}$ as a subsequence
 $\fP=(d_{k+1},d_{k+2},\dots,d_{k+l})$, where $0\le k< rs,\ 1\le l\le rs$, such that $d_i\ge0$ for all $i=k+1,\dots,k+l$, but either
 $l=rs$ or both $d_k<0$ and $d_{k+l+1}<0$ (again we set $d_{rs+j}=d_j$). Set $\th(\fP)=l$ if either $l=rs$ or $\fP=(0,0,\dots,0)$
 and $\th(\fP)=l+1$ otherwise; $\th(\fD)=\sum_{\fP}\th(\fP)$.
 At last, set $\de(\fD,\la)=1$ if $\fD={\bf0}=(0,0,\dots,0),\, \la=1$, and $\de(\fD,\la)=0$ otherwise. Then 
 \begin{align*} 
  \nH^0(E,\kV(\fD,m,\la))&=m\left(\sum_{i=1}^{rs}(d_i+1)^+\! -\th(\fD)\right)+\de(\fD,\la),\\
  \nH^1(E,\kV(\fD,m,\la))&=m\left(\sum_{i=1}^{rs}(d_i+1)^-\!+rs -\th(\fD)\right)+\de(\fD,\la).
 \end{align*} 
 Moreover, $\kV(\fD,m,\la)$ is generically spanned \iff either $\fD=\bf0$, $m=1$ and $\la=1$ (i.e.
 $\kV(\fD,m,\la)\iso\kO_E$) or the following conditions hold%
\footnote{I must note a mistake in the preprint \cite{dgk}, where we claimed that $\fD>\bf0$ is enough
 for $\kV(\fD,m,\la)$ to satisfy Kahn's conditions. It has been improved in the final version.
 We are thankful to Igor Burban who has noticed this mistake.}:
 \begin{enumerate}
\item  $\fD>\bf0$  (it means that all $d_i\ge0$ and at least one inequality is strict).
 \item 
 If $\fD'$ is a shift of $\fD$, i.e. $\fD'=(d_{k+1},\dots,d_{rs},d_1,\dots,d_k)$ for some $k$,
  $\fD'$ contains no subsequence $(0,1,1,\dots,1,0)$, in particular $(0,0)$.
 \item
  No shift of $\fD$ is of the form $(0,1,1,\dots,1)$.
\end{enumerate}
 We call a sequence $\fD$ satisfying the conditions 1--3 a \emph{suitable sequence}. The inequality $\fD>\bf0$ implies
 that $\rH^1(E,\kV(\fD,m,\la))=0$, so the vector  bundles $\kV(\fD,m,\la)$ with suitable $\fD$ satisfy Kahn's conditions.
 Moreover, in this case $Q=\rH^1(E,\kV_\fD(Z))\iso\rH^1(E,\kV_{\fD-\fB^r})$ is
 flat over $\bR=\Mk[x,x^{-1}]$ if $\fD\ne\fB$. If $\fD=\fB$, $Q$ is no more flat, but $Q[(x-1)^{-1}]$ is flat over 
 $\bR'=\Mk[x,x^{-1},(x-1)^{-1}]$. Set $\kV'_\fB=\kV_\fB\*_\bR\bR'$ and 
 $ n_\fD=\sum_{i=1}^{rs}(d_i-b_i+1)^+\!-\th(\fD-\fB^r)$.
 Note that $n_\fD=0$ \iff every positive part of $\fD-\fB^r$ contains $1$ at most once,
 all other entries of it being $0$ (for instance $\fD=\fB$).
 Define, for each suitable sequence $\fD$, a family $\kE_\fD$ of vector bundles over $E$ based on $\bR$ if $\fD\ne\fB$ and 
 on $\bR'$ if $\fD=\fB$. Namely,
 \begin{align*}
 \kE_\fD&=\kV_\fD\+n_\fD\kO_E\*\bR, \text{ if } \fD\ne\fB;\\
 \kE_\fB&=\kV'_\fB.
 \end{align*} 
  These families satisfy the conditions of Corollary \ref{c73}, hence can be lifted 
 to full families $\kF_\fD$ of vector bundles over $X$, which give rise to families $\kM_\fD=\pi_*\kF_\fD$ of \CM $\bA$-modules.
 Now Theorems \ref{t41} and \ref{t72} directly imply the following description
 of \CM $\bA$-modules.

 \begin{theorem}\label{t77}
  If $\bA$ is a cusp singularity, all indecomposable $\bA$-modules are:
 \begin{itemize}
\item 
  $\bA$, of rank 1;
 \item
  $\kM_\fD(m,\la)=\kM_\fD\*_\bR\bR/(x-\la)^m$, of rank $m(r+n_\fD)$, where $\fD$
 is a suitable sequence and $\la\ne1$ if $\fD=\fB$;
 \item
  $\kM_\fB(m,1)$, of rank $m+1$, where $R_E\kM_\fB(m,1)\iso\kV(\fB,m,1)\+\kO_E$.
\end{itemize}
 In particular $\bA$ is \CM tame unbounded.
 \end{theorem} 
 \end{proof} 
 
 \section{\CM types: Q-elliptic case}
 \label{sec12}

 Using Corollary \ref{c54} we can extend the results of Section \ref{sec11} to a wider class of surface singularities.

 \begin{defin}\label{d121}
  \begin{enumerate}
\item 
 A surface singularity $\bA$ will be called \emph{Q-Gorenstein}, if the order $g$ of its canonical divisor $K$ in the Picard group
 is finite and prime to $\chr\Mk$.

 Just as in the proof of Theorem \ref{t54}, one can construct the \emph{Gorenstein cover} $\bB=\bop_{i=1}^gJ^{-i}t^i$, where
 $J$ is an ideal of class $K$ and $t^n=\al^{-1}$ such that $J^{-n}=\al\bA$. 
 \item
  A Q-Gorenstein singularity $\bA$ is called
 \begin{itemize}
\item 
 \emph{Q-elliptic} if its Gorenstein cover is minimally elliptic;
 \item
  \emph{simple Q-elliptic} if its Gorenstein cover is simple elliptic;
 \item
  \emph{Q-cusp} if its Gorenstein cover is a cusp singularity.
\end{itemize}  
\end{enumerate}
 \end{defin} 

 Note that if $\bA$ is Q-Gorenstein and $\bB$ is its Gorenstein cover, the cyclic group of
 order $g$ acts on $\bB$ so that $\bA=\bB^G$ and the extension $\bA\sbe\bB$ is unramified in codimension 1 (cf. the proof of
 Theorem \ref{t54}). Therefore Theorem \ref{t75} and Corollary \ref{c54} immediately imply

 \begin{corol}\label{122}
  A Q-elliptic singularity is
 \begin{itemize}
\item 
 \CM \emph{tame bounded} if it is simple Q-elliptic;
 \item
  \CM \emph{tame unbounded} if it a Q-cusp singularity;
 \item
  \CM \emph{wild} otherwise.
\end{itemize}
 \end{corol} 

 \begin{remk}
  If $\Mk=0$, simple Q-elliptic and Q-cusp singularities arise as the so called \emph{log-canonical} singularities 
 \cite{kaw}. Recall that a normal surface singularity with the minimal resolution
 $\pi:X\to S$ is called \emph{log-canonical} if $K_X=\pi^*K_S+\sum_{i=1}^sa_iE_i$ with $a_i\ge-1$, where $\pi^*K_S$
 denotes the \emph{numerical pullback} of $K_S$. The latter means that $\pi^*K_S=K'_S+\sum_{i=1}^sk_iE_i$, where
 $K'_S$ is the strict transform of $K_S$ (cf. \cite[Section II.7]{ha}) and $k_i\in\mathbb Q$ are so chosen that $(\pi^*K_S,E_i)=0$
 for all $i$. In other words, $K_X$ coincides with $K_S$ outside $E$, while $(K_X.E_j)=\sum_{i=1}^sa_i(E_i.E_j)$ for all $j$.
 It is proved in \cite[Theorem 9.6]{kaw} that, if $\chr\Mk=0$, a log-canonical singularity is always Q-Gorenstein;
 moreover, it is either quotient, or simple Q-elliptic, or Q-cusp.
 \end{remk} 
 
 We shall describe \CM modules over Q-cusp singularities. Namely, let $\bB$ be a
 cusp singularity,  $T=\spec\bB$, $\phi:Y\to T$ be the minimal resolution of $T$, and $F=\bup_{i=1}^tF_i$ be the exceptional
 curve of $\phi$, where $F_i$ are the irreducible components of $F$ so arranged that $F_i\cap F_{i+1}=\0$ (as usually we
 set $F_{t+i}=F_i$ and use analogous identification everywhere). Suppose that $\chr\Mk\ne2$ and the group $G=\set{1,\si}$
 of order 2 acts on $\bB$ so that the lifting of this action onto $Y$ is free outside $F$ and reverse the orientation of $F$
 (Q-cusp singularities fit this situation: it can be shown just as in \cite[Theorem 9.6]{kaw}). Then $\si$ induces a reflection
 of $\De(F)$. If we choose $F_1$ to be its fixed component, $\si$ maps $F_i$ onto $F_{t+2-i}$ with the coordinate
 transformation $x\mps1/x$. It induces the action of $\si$ on vector bundles over $F$ such that
 $\kV(\fD,m,\la)^\si\iso\kV(\fD^\si,m,1/\la)$, where $\row d{rt}^\si=(d_1,d_{rt},d_{rt-1},\dots,d_2)$.Therefore
 its action on \CM $\bB$ modules is (using notation of Theorem \ref{t77}):
 $$
   \kM_\fD(m,\la)^\si\iso\kM_{\fD^\si}(m,1/\la),\ \bB^\si\iso\bB\ \text{ and }\ M_n^\si\iso M_n
 $$ 
 (since necessarily $\fB^\si=\fB$). But it is known from generalities about group actions that the restriction onto $\bA$
 of  an indecomposable $\bB$-module $M$ is
 indecomposable if $M^\si\not\iso M$ and decomposes $M=M'\+M''$, where $M',M''$ are indecomposable
 and non-isomorphic, if $M^\si\iso M$. Moreover, if $N$ is another module, $N\not\iso M$ and $N\not\iso M^\si$, then $M^\si$
 and $N^\si$ have no common direct summands.  Denote by $\kN_\fD$ the following families of \CM $\bA$-modules:
 \begin{itemize}
\item 
 the restriction of $\kM_\fD$ onto $\bA$ if $\fD^\si$ does not coincide with any $t$-shift of $\fD$;
 \item
  the restriction of $\kM_\fD[(x^2-1)^{-1}]$ onto $\bA$ if $\fD^\si$ is a $t$-shift of $\fD$,
\end{itemize}
 where $\bR=\Mk[x,x^{-1}]$. For any sequence $\fD$ such that $\fD^\si$ coincides with a $t$-shift of $\fD$ denote by
 $N_\fD'(m,\pm1)$ and $N_\fD''(m,\pm1)$ the indecomposable direct summands of $\kM_\fD(m,\pm1)$. Note that 
 always $\fB^\si=\fB$. The previous observations together with Theorem \ref{t77} imply 

 \begin{theorem}\label{t123}
  In the above situation all indecomposable \CM $\bA$-modules are:
 \begin{itemize}
\item
  $\kN_\fD(m,\la)=\kN_\fD\*_\bR\bR/(x-\la)^m$, where $\la\ne0$, and $\la\ne\pm1$ if $\,\fD^\si$ is a $t$-shift of $\,\fD$;
 \item
  $\kN_\fD'(m,\pm1)$ and\, $\kN_\fD''(m,\pm1)$ for such $\fD$ that $\,\fD^\si$ is a $t$-shift of $\,\fD$;
 \item
  $\bA$ and\, $\bB^-=\setsuch{b\in\bB}{\si(b)=-b}$.
\end{itemize}
\em Here $\fD$ always denotes a suitable sequence.
 \end{theorem} 

 \section{Application to hypersurfaces and curves.}
 \label{sec13}
 
 We can apply the results of Section \ref{sec11} to \emph{hypersurface singularities}, i.e. rings of the shape
 $\bA=\Mk[[\lst xn]]/(f)$, using the results of Kn\"orrer \cite{kn} (see also \cite{yo}) on the relations
 between \CM modules over $\bA$ and over its \emph{suspension} $\bA\vs=\Mk[[\lst x,z]]/(f+z^2)$.
 (In this section we suppose that $\chr\Mk\ne2$.)
 Namely, for every $\bA$-module $M$ denote by $\syz M$ its first syzygy as of $\bA$-module
 and by $\Om M$ its first syzygy as of $\bA\vs$-module. For every $\bA\vs$-module
 $N$ denote by $\res N$ the $\bA$-module $N/zN$. These operations map \CM modules to \CM ones,
 so we consider them as functors between categories of \CM $\bA$- and $\bA\vs$-modules. As $\bA$
 is Gorenstein, $\syz$ can be considered as an automorphism of the \emph{stable category}
 $\ul{\Cm}(\bA)$, which is the factor of $\Cm(\bA)$ modulo free modules.
 We also denote by $\si$ the automorphism of $\bA\vs$ mapping $z$ to $-z$ and
 leaving all $x_i$ fixed, and by $N^\si$ the $\bA\vs$-module obtained from $N$ by twisting with $\si$.
 Then the results of \cite{kn} can be formulated as follows.

 \begin{theorem}\label{131}
  \begin{enumerate}
 \item
  Let $M\not\iso\bA$ be an indecomposable \CM $\bA$-module. 
  \begin{enumerate}
\item 
  If $M\not\iso \syz M$, the module $\Om M$ is indecomposable.
 \item
  If $M\iso\syz M$, $\Om M\iso \Om_1M\+\Om_2M$, where $\Om_1M$ and $\Om_2M$ are
 indecomposable and non-isomorphic.
\end{enumerate}
 Every indecomposable \CM $\bA\vs$-module is isomorphic to one of those described in items {\rm (a),(b)}.
\item 
 Let $N\not\iso\bA\vs$  be an indecomposable \CM $\bA\vs$-module. 
  \begin{enumerate}
\item 
 If $N\not\iso N^\si$, the module $\res N$ is indecomposable.
 \item
  If $N\iso N^\si$, $\res N=\res_1N\+\res_2N$, where $\res_1N$ and $\res_2N$ are indecomposable
 and non-isomorphic.
\end{enumerate}
 Every indecomposable \CM $\bA$-module is isomorphic to one of those described in items {\rm (a),(b)}.
\end{enumerate}
 \end{theorem} 
 
 We apply these results to singularities of type $\rT_{pqr}$. Namely, denote by $\bT_{pqr}$ the factor
 $$
  \Mk[[x,y,z]]/(x^p+y^q+z^r+\la xyz),\quad \text{where }\ 1/p+1/q+1/r\le1\ \text{ and }\ \la\ne0
 $$
 (we may suppose that $p\ge q\ge r$).
 Moreover, we demand this singularity to be \emph{isolated}, which imposes restrictions on $\la$ in the \emph{quasi-homogeneous}
 cases, when $(p,q,r)\in\set{(2,3,6),(2,4,4),(3,3,3)}$. Note that in all other cases the isomorphism class of $\bT_{pqr}$ does not
 depend on $\la$. A hypersurface singularity that is an (iterated) suspension of $\bT_{pqr}$,
 i.e. $$\Mk[[x,y,z,\lst tm]]/\left(x^p+y^q+z^r+\la xyz+\sum_{i=1}^mt_i^2\right)$$
 is called \emph{a singularity of type} $\rT_{pqr}$. It is known \cite{lau} that the surface singularities $\bT_{pqr}$ are simple elliptic
 in quasi-homogeneous case and cusp singularities in all other cases. Therefore they are tame, so Theorem \ref{131} implies

 \begin{corol}\label{132}
  All hypersurface singularities of type $T_{pqr}$ are \CM tame.
 \end{corol} 
 Unfortunately, we do not have precise formulae for syzygies of $\bT_{pqr}$-modules, so we cannot give an explicit description of
 \CM modules over their suspensions. 

 We can also use Kn\"orrer's correspondence to obtain a description of \CM modules over curve singularities of types $\rT_{pq}$,
 i.e. the factors
 $$
  \bT_{pq}=\Mk[[x,y]]/(x^p+y^q+\la x^2y^2),\quad \text{where }\ 1/p+1/q\le1/2\ \text{ and }\ \la\ne0
 $$
 (again in quasi-homogeneous cases, when $(p,q)\in\set{(3,6),(4,4)}$, some conditions must be imposed on $\la$ in order this
 singularity to be isolated). To do it one only has to note that the suspension of $\bT_{pq}$ is isomorphic to a 
 surface singularities $\bT_{pq2}$ (change $z$ to $z+\la xy/2$). Moreover, in this case one can explicitly calculate the
 action of $\si$ on the minimal resolution of $\bT_{pq2}$, thus get an explicit list of indecomposable $\bT_{pq}$-modules
 (cf. \cite{dgk} for details). It accomplishes the study of \CM modules over tame curve singularities, filling the flaw in \cite{dg1},
 where no explicit description of $\bT_{pq}$-modules was obtained (their tameness was established using deformation theory). 

 \section{Some conjectures and remarks}
 \label{sec14}

 We end up with some conjectures and remarks (cf. \cite{dgk}).

 \begin{conj}\label{141}\brn
  In the following cases the ring $\bA$ is \CM wild:
 \begin{enumerate}
\item 
 $\bA$ is a surface singularity that is neither quotient, nor simple Q-elliptic, nor Q-cusp;
 \item
  $\bA$ is a hypersurface singularity that is neither simple (i.e. of types {\rm A-D-E} \cite{avg})
 nor of type $\rT_{pqr}$;
 \item
  $\bA$ is non-isolated and the dimension of its singular locus is greater than 1 (i.e. $\bA_\gP$ is not regular for some
 prime ideal $\gP$ of depth 2).
\end{enumerate}
 \end{conj} 

 If this conjecture is true, we shall have a complete description of \CM types of isolated surface and hypersurface singularities.
 The result is given in Table 1 (the conjectured cases marked with `?'). 
 Unfortunately, we have now no further conjectures, not even examples, for non-isolated singularities with 1-dimensional
 singular locus. Probably, very few of them can be \CM tame. 

\begin{table}
 \begin{center}
 Table 1.\\\medskip \sf Cohen--Macaulay types of singularities
 \end{center}
\bigskip 
\begin{center}
\begin{tabular}{|c||c|c|c|}
\hline &&&\\
 CM type & curves & surfaces & hypersurfaces \\
&&&\\ \hline\hline &&&\\
 finite & dominate  & quotient & simple \\
 & A-D-E && (A-D-E) \\
&&&\\ \hline &&&\\
 tame & dominate  & simple &  of type $\rT_{pqr}$ with\\
 bounded & $\bT_{pq}$ with    & Q-elliptic &  $1/p+1/q+1/r=1$ \\
 &  $1/p+1/q=1/2$    & (only\,?)  &(only\,?)  \\
&&&\\ \hline &&&\\
  tame & dominate  & Q-cusp & of type $\rT_{pqr}$ with \\
 unbounded & $\bT_{pq}$ with    &(only\,?)  & $1/p+1/q+1/r<1$  \\
& $1/p+1/q<1/2$      &  &(only\,?) \\
&&&\\ \hline &&&\\
 wild & all other & {}\ all other\,?\ {}& all other\,? \\
&&&\\ \hline 
\end{tabular}
\end{center}
\end{table}

\begin{remk}\label{142}
 All known examples of \CM tame unbounded singularities, in particular those from Table 1,
 are actually of exponential growth.  It seems very plausible that it is always so. Nevertheless,
 just as in the case of finite dimensional algebras, it can only be shown \emph{a posteriori},
 when one has a description of modules. We do not see any ``natural'' way to prove this conjecture
 without such calculations. That is why we prefer to say `bounded' and  `unbounded' instead of more usual
 `of polynomial  growth' and `of exponential growth.'
\end{remk}

\medskip
\begin{remk}
\label{143}
 In the complex analytic case, Artin's Approximation Theorem \cite{art} implies
 that the list of \CM modules remains the same if
 $\bA$ denotes the ring of germs of analytic functions on a simple elliptic or cusp 
 singularity. The lifting of families is more cumbersome. It is always possible if the base
 is a finite dimensional algebra. If it is an algebraic variety $T$,
 we can only claim that for each point $t\in T$ a lifting is possible over a neighbourhood $U$ of
 $t$ in $T$.  It gives a lifting of an appropriate family to an etale covering $\ti T$ of $T$. If $T$ is a smooth
 curve or surface, so is $\ti T$, therefore the results on tameness and wildness from
 Sections~\ref{sec11},\ref{sec12} remain valid. On the other hand, in the case of cusps it seems credible
 that the families $\kE_\fD$ from the proof of Theorem~\ref{t75}
 can actually be lifted over $T$, just as in \cite{kahn} for simple elliptic case.
\end{remk}

\end{document}